\documentclass[11pt, letterpaper]{amsart}
\usepackage{amsmath, amssymb, amsthm,fullpage}
\usepackage[all]{xypic}
\newtheorem{thm}[equation]{Theorem}
\newtheorem{prop}[equation]{Proposition}
\newtheorem{cor}[equation]{Corollary}
\newtheorem{lemma}[equation]{Lemma}

\newtheorem{ques}[equation]{Question}

\newtheorem{claim}[equation]{Claim}
\theoremstyle{definition}
\newtheorem{defn}[equation]{Definition}

\theoremstyle{remark}
\newtheorem{exam}[equation]{Example}

\newtheorem{ntn}[equation]{Notation}

\newtheorem{rem}[equation]{Remark}

\makeatletter
\renewcommand{\subsection}{\@startsection{subsection}{2}{0pt}{-3ex
plus -1ex minus -0.2ex}{-2mm plus -0pt minus
-2pt}{\normalfont\bfseries}} 
\renewcommand{\subsubsection}{\@startsection{subsubsection}{2}{0pt}{-3ex
plus -1ex minus -0.2ex}{-2mm plus -0pt minus
-2pt}{\normalfont\bfseries}} \makeatother

\numberwithin{equation}{section}
\numberwithin{equation}{subsection}

\newdir^{ (}{{}*!/-5pt/@^{(}}%  % A new tip style for right hook arrows.
\newcommand{\tors}{\operatorname{Torsion}}
\newcommand{\Ind}{\operatorname{Ind}}

\newcommand{\norm}[1]{\left\| #1 \right\|}
\newcommand{\Ker}{\operatorname{Ker}}
\newcommand{\imm}{\operatorname{im}}
\newcommand{\coker}{\operatorname{coker}}

\newcommand{\LL}{\mathcal{L}}

\newcommand{\Tor}{\operatorname{Tor}}
\newcommand{\Ext}{\operatorname{Ext}}
\newcommand{\uAmd}{\underline{A-mod}}

\newcommand{\uAmk}{\underline{A-mod}_{\k}}

\newcommand{\Stab}{\operatorname{Stab}}
\newcommand{\uHom}{\operatorname{\underline{Hom}}}
\newcommand{\uExt}{\operatorname{\underline{Ext}}}
\newcommand{\uTor}{\operatorname{\underline{Tor}}}

\newcommand{\uRHom}{\operatorname{\underline{RHom}}}

\newcommand{\uo}{\operatorname{\underline{\otimes}}}
\newcommand{\uHH}{\operatorname{\underline{HH}}}
\newcommand{\HH}{\operatorname{HH}}

\newcommand{\Dk}{\mathbb{D}_{\k}}

\newcommand{\RHom}{\operatorname{RHom}}

\newcommand{\tra}{\mathop{\rightarrow}}

\newcommand{\qiso}{{\;\stackrel{q.i.}{\stackrel{_\sim}{\to}}\;}}
\newcommand{\iso}{{\;\stackrel{_\sim}{\to}\;}}
\newcommand{\liso}{{\;\stackrel{_\sim}{\leftarrow}\;}}

\newcommand{\beq}{\begin{equation}\label}
\newcommand{\D}{\mathbb{D}}

\newcommand{\eeq}{\end{equation}}

\newcommand{\into}{\hookrightarrow}

\newcommand{\onto}{\twoheadrightarrow}

\newcommand{\ds}{\displaystyle}
\newcommand{\Hom}{\operatorname{Hom}}

\newcommand{\op}{\text{op}}

\newcommand{\V}{\mathcal{V}}
\newcommand{\W}{\mathcal{W}}
\newcommand{\rptop}{r_{p,\mathrm{top}}}

\def\R{\mathbb{R}}
\def\C{\mathbb{C}}
\def\k{\mathbf{k}}

\def\o{\otimes}

\def\End{\mathrm{End}}
\def\dq{{\overline{Q}}}

\def\Id{\mathrm{Id}}

\def\N{\mathbb{Z}_{\geq 0}}
\def\Q{\mathbb{Q}}
\def\F{\mathbb{F}}

\def\dim{\operatorname{dim}}

\def\Z{{\mathbb Z}}
\def\tr{{\operatorname{tr}}}

\def\1{\mathbf{1}}

\def\ldp{\mathopen{(\!(}} \def\rdp{\mathclose{)\!)}}

\begin{document}
\title{Calabi-Yau Frobenius algebras}
\author{Ching-Hwa Eu and Travis Schedler} 
\maketitle

\begin{abstract}
  We define Calabi-Yau Frobenius algebras over arbitrary
  base commutative rings.  We define a Hochschild analogue of Tate
  cohomology, and show that this ``stable Hochschild cohomology'' of
  periodic CY Frobenius algebras has
  a Batalin-Vilkovisky and Frobenius algebra structure.
   Such algebras include (centrally extended) preprojective
  algebras of (generalized) Dynkin quivers, and group algebras of
  classical periodic groups.  We use this theory to compute (for the
  first time) the Hochschild cohomology of many algebras related to
  quivers, and to simplify the description of known results.
  Furthermore, we compute the maps on cohomology from extended Dynkin
  preprojective algebras to the Dynkin ones, which relates our CY
  property (for Frobenius algebras) to that of Ginzburg (for algebras
  of finite Hochschild dimension).
\end{abstract}
\tableofcontents

\section{Introduction}
Frobenius algebras have wide-ranging applications to geometry, e.g.,
to TQFTs in \cite{Moo,Laz,KC1}, and are closely related to the string
topology Batalin-Vilkovisky (BV) algebra \cite{CS,CS2}, as described
in e.g.~\cite{KC1,HL,CV,TZ}. 

One motivation for this work is to explain the following algebraic
phenomenon: the Hochschild cohomology of many interesting Frobenius
algebras has a BV structure.  For example, this is true for symmetric
Frobenius algebras, and for preprojective algebras of Dynkin quivers
(which are \textit{not} symmetric).  

To explain this, we define \textit{Calabi-Yau (CY) Frobenius
  algebras},
whose Hochschild cohomology has not only the usual Gerstenhaber
algebra structure, but a BV structure (at least when the algebra is
also periodic, as we will explain).  The CY Frobenius property is
similar to the CY property of \cite{Gcy}, but not the same: any
Frobenius algebra over $\C$ which is CY in the sense of \cite{Gcy}
must be a direct sum of matrix algebras.

The prototypical example of a CY Frobenius algebra $A$ \textit{of dimension $m$} ($CY(m)$) over a field
$\k$ is one that has a resolution (for $A^e := A \o_\k A^{\op}$)
\begin{equation}
A^\vee \into P_m \rightarrow P_{m-1} \rightarrow \cdots \rightarrow P_0 \onto A, \quad A^\vee := \Hom_{A^e}(A, A^e),
\end{equation}
with each $P_i$ a projective $A$-bimodule (the sequence must be
exact).  Our definition is a weakened version of this:
$A^\vee \simeq \Omega^{m+1} A$ in the stable $A$-bimodule category (this
notion is recalled in Appendix \ref{gfs}). 

Any symmetric Frobenius algebra has $A \cong A^\vee$, and so we say it
is CY Frobenius of dimension $-1$. Our methods give a very short,
simple proof of the fact that the Hochschild cohomology of such
algebras is BV (\S \ref{symfrobs}), for which many (generally
more complicated) proofs are given in
e.g.~\cite{KC1,T,MBV,TZ,K04,K06}.

In \cite{ES,ES2,EE2}, it was also noticed that the Hochschild
cohomology of preprojective algebras of Dynkin quivers has a
self-duality property: one has $HH^{i}(A) \cong HH^{5+6j-i}(A)^*$ for
such algebras, with $i, 5+6j-i \geq 1$, and moreover one has the
periodicity $HH^i(A) \cong HH^{i+6j}$ when $i, 6j+i \geq 1$ (here, $j$
is an integer).

We explain this by introducing \textit{stable Hochschild
  cohomology}
$\uHH^\bullet(A)$, of Frobenius algebras $A$. This is a Hochschild
analogue of Tate cohomology, which coincides with the usual Hochschild
cohomology in positive degrees, and is defined using the stable module
category by $\uHH^\bullet(A) = \uHom(\Omega^\bullet A, A)$ (in the
latter form, this was studied in various papers, e.g., \cite{ESk}). We
show that this is a $\Z$-graded ring, and prove it is
graded-commutative.  To our knowledge, this is the first time
$\uHom(\Omega^\bullet A, A)$ has been studied in this way.  We also
define the notion \textit{stable Hochschild homology}
$\uHH_{\bullet}(A)$, which coincides with the usual notion in positive
degrees, and prove that $\uHH_{\bullet}(A)$ is a graded module over
$\uHH^\bullet(A)$ using a natural \textit{contraction} action, which
generalizes contraction in the $\Z_{\geq 0}$-graded case.  These
results are in Theorem \ref{asfhh}.

We then show that, in general,
$\uHH_{\bullet}(A) \cong \uHH_{-1-\bullet}(A)^*$ (which makes the
contraction operations \textit{graded self-adjoint}), and that the
$CY(m)$-Frobenius property produces dualities
$\uHH_\bullet(A) \cong \uHH^{m-\bullet}(A)$.  Put together, we
obtain a \textit{graded Frobenius algebra structure} on
  $\uHH^\bullet(A)$
(Theorem \ref{cyfhhthm}), explaining the aforementioned results of
\cite{ES,ES2,EE2}. In particular, $\uHH^\bullet(A) \cong \uHH^{2m+1-\bullet}(A)^*$, as modules over $\uHH^0(A)$; in the preprojective algebra cases, $m=2$,
which explains the aforementioned duality.

Call a (Frobenius) algebra \textit{periodic} if it has a periodic
$A$-bimodule resolution (or, more generally, $A \simeq \Omega^n A$ in
the stable bimodule category, for some $n$).  We prove that \textit{the stable
  Hochschild cohomology of periodic CY Frobenius algebras has a BV
  structure}
(Theorem \ref{pcyfralg}).  More generally, for any periodic (not
necessarily CY) Frobenius algebra, the structure of calculus
\cite{TT,DGT} on the pair $(\HH^\bullet, \HH_\bullet)$ extends to the
$\Z$-graded setting $(\uHH^\bullet, \uHH_\bullet)$ (Theorem
\ref{perfhhthm}). Moreover, we show that, in the (centrally extended)
preprojective cases, the BV structure and Frobenius algebra structure
are compatible: the BV differential is graded selfadjoint---we
call such an algebra a \textit{BV Frobenius algebra}.

Additionally, all of the above work is done in the context not only of
Frobenius algebras over a field, but over an arbitrary base
commutative ring.  Precisely, we use the notion of \textit{Frobenius
  extensions of the first kind}
\cite{NT2,Ka2}, which says that the algebra is projective over a base
commutative ring $\k$ and that the duality is nondegenerate over this
base.  (Perhaps this could be generalized further to an
$A_\infty$-Frobenius property, but we do not do this here.)

New computational results (extending \cite{ES,ES2,EE2,Eu2,Eu3})
concerning preprojective algebras include:
\begin{itemize}
\item The computation of Hochschild (co)homology of preprojective
  algebras of Dynkin type over $\k=\Z$ (in types $D,E$, the results
  \cite{EE2,Eu2,Eu3} are over characteristic-zero fields), and the
  proof that $\uHH^\bullet$ is \textit{BV Frobenius} over any field
  (Theorem \ref{prepthm});
\item The computation of cup product and BV structure on Hochschild
  (co)homology of \textit{centrally extended preprojective algebras}
  \cite{ER} (the groups were computed over $\C$ in \cite{Eu}), which
  we show is \textit{BV Frobenius} (Theorem \ref{cprepthm});
\item The computation of the induced maps from the Hochschild
  (co)homology in the \textbf{extended} Dynkin case to the Dynkin case
  by cutting off the extending vertex (Theorem \ref{eddcompthm}). This
  explains the structure in the Dynkin case and elucidates the
  relationship between the usual Calabi-Yau and Calabi-Yau Frobenius
  properties (which is analogous to Euclidean vs. spherical geometry);
\end{itemize}
We also explain and simplify the cited known results.

Our other main example is the case of group algebras of finite groups,
which are automatically Calabi-Yau Frobenius (since they are
symmetric).  We are interested in when these are periodic (and hence
$\uHH^\bullet$ is BV and Frobenius, by Theorem \ref{pcyfralg}). We
show (Theorem \ref{gpalgthm}) that the periodic (CY) Frobenius
algebras are just the classical periodic algebras, i.e., those whose
group cohomology is periodic (using classical results).  In the
appendix, we also give an elementary topological proof that groups
that act freely and simplicially on a sphere (such as finite subgroups
of $SO(n)$) are periodic (CY) Frobenius.

\subsection{Acknowledgements}
We thank P. Etingof and V. Ginzburg for useful discussions
and advice.  The second author was partially supported by an NSF GRF.

\subsection{Definitions and Notation}
Here we recall some standard definitions and state the notation we will
use throughout.

All complexes will be assumed to have decreasing degree (i.e., they
are chain complexes), unless otherwise specified.

Let us fix, once and for all, a commutative ring $\k$.  
  When we say ``algebra over $\k$'', we mean an algebra $A$ over
  $\k \rightarrow A$ such that the image of $\k$ is central in $A$.  
  Bimodules over an algebra $A$ over $\k$ will be assumed to be symmetric as
  $\k$-bimodules (i.e., $A$-bimodules mean $A^e := A \o_\k A^{\op}$-modules).
\begin{ntn} The category $A-mod$ means \textit{finitely-generated} $A$-modules,
for any ring $A$.  The category $A-mod_{\k}$ means finitely-generated $A$-modules \textit{which are finitely-generated projective as $\k$-modules}.
\end{ntn}
\begin{ntn}
We will abbreviate ``finitely-generated'' as ``fg.''
\end{ntn}

By a ``Frobenius algebra over a commutative ring $\k$'', we will mean what
is also known as a ``Frobenius extension of the first kind'' in the literature
\cite{NT2,Ka2}, namely:
\begin{defn}
  A Frobenius algebra $A$ over a commutative ring $\k$ is a $\k$-algebra
  which is a fg projective $\k$-module, and which
  is equipped with a nondegenerate invariant inner product $(\,, )$,
  i.e.:
\begin{equation}
(ab, c) = (a, bc), \forall a,b,c \in A; \quad (-, a): A \rightarrow \Hom_\k(A, \k) \text{ is an isomorphism of $\k$-modules}.
\end{equation}
\end{defn}
\begin{exam}
  Any group algebra $\k[G]$ is Frobenius if $G$ is finite, using the
  pairing $(g, h) = \delta_{g,h^{-1}}$ for $g, h \in G$. In fact, this algebra
  is symmetric (meaning $A \cong A^*$ as $A$-bimodules).
\end{exam}
\begin{exam}
  The preprojective algebra $\Pi_Q$ is known to be Frobenius if $Q$ is
  Dynkin (cf.~e.g., \cite{ES2, ER}). It is not difficult to see (e.g., through explicit bases as in \cite{Eu2})
  that these are in fact Frobenius over $\Z$.
\end{exam}
\begin{exam}
  Centrally extended preprojective algebras of Dynkin quivers were
  defined and proved to be Frobenius in \cite{ER} (working over $\C$).
\end{exam}
\begin{exam}
  For any two Frobenius algebras $A, B$ over $\k$, the algebra
  $A \o_\k B$ is also Frobenius. In particular, so is $A^e$.
\end{exam}
We will mainly be interested in the cases where $\k$ is a field or
$\k = \Z$, and $A$ is a fg free $\k$-module. (These
are generally known as ``free Frobenius extensions of the first
kind,'' cf.~\cite{NT2}.)

We refer to Appendix \ref{gfs} for some general results about
Frobenius algebras and the stable module category over arbitrary
commutative rings, which are direct generalizations of standard
results in the case where $\k$ is a field.   In particular, the results
there justify the following definitions:
\begin{defn}
Let the projectively stable module category $\uAmd$ be the category whose objects
are fg $A$-modules,
and whose morphisms $\uHom_{A}$ are given by
\begin{multline}
\uHom_{A}(M, N) := \{f \in \Hom_A(M,N)
\} \\ /\{\text{morphisms that factor
through a projective $A$-module}\}.
\end{multline} 
Let $\uAmk \subset \uAmd$ be the full subcategory of modules which are
projective as $\k$-modules.
\end{defn}

\begin{defn}
For any algebra $A$ over $\k$, define the functor $^\vee: A-mod \rightarrow A^{\op}-mod$ by
\begin{equation} \label{veedfn}
M^\vee := \Hom_A(M, A),
\end{equation}
with the natural induced maps on morphisms.
\end{defn}

\begin{defn}
For any Frobenius algebra $A$, let $\eta: A \iso A$ be the \textbf{Nakayama automorphism}
defined by
\begin{equation}
(a, b) = (\eta^{-1}(b), a), \quad \forall a, b \in A.
\end{equation}
\end{defn}
\begin{defn} For any ($\k$-linear) automorphism $\phi: A \iso A$, and
  any $A$-module $M$, let $_\phi M$ denote $M$ with the twisted action
  given by precomposing $A \rightarrow \End_\k(M)$ by $\phi$.
  Similarly, for any bimodule $N$, $_\phi N_\psi$ denotes twisting the
  left action by $\phi$ and the right action by $\psi$.  If either of
  $\phi, \psi$ is the identity, we may omit it from the notation.
\end{defn}
We have $A_\eta \cong A^*$ as $A$-bimodules (conversely, such an
isomorphism is equivalent to $(-,-)$ with automorphism $\eta$).

\begin{defn} Let $A$ be a Frobenius algebra over $\k$ and $M, N$
  fg $A$-modules which are projective as $\k$-modules.
  For any integer $i \in \Z$, let us denote
\begin{equation}
\uExt_A^i(M, N) := \uHom_A(\Omega^i M, N).
\end{equation}
\end{defn}
\begin{rem} \label{taterem} When $G$ is a finite group, and $\k$ is a
  field (of any characteristic), the cohomology groups
  $\uExt^i_{\k[G]}(\k, M)$ are the Tate cohomology groups over
  $\k$ with coefficients in $M$.  Indeed, we may compute
  $\uExt^i_{\k[G]}(\k,M)$ by the complex $\Hom(P_\bullet, M)$ where
  $P_\bullet$ is any two-sided projective resolution of $\k$ (i.e., an
  exact (unbounded) complex of projectives such that the cokernel of
  $P_{1} \rightarrow P_0$ is $\k$).  This is one of the standard
  definitions of Tate cohomology (cf.~e.g.~\cite{AM}, Definition 7.1).
\end{rem}
\begin{rem}
  By the same token, it makes sense to define the Tate cohomology of
  any Hopf algebra $H$ which is Frobenius over $\k$ by
  $\uExt^i_{H}(\k, M)$, where $\k$ is the augmentation module. Note
  that, if $\k$ is a PID, then a Hopf algebra $H$ over $\k$ is
  automatically Frobenius if it is fg projective as a
  $\k$-module \cite{LS}.
\end{rem}

Next, we recall the definition of a Gerstenhaber algebra.
\begin{defn}
A Gerstenhaber algebra $(\mathcal{V}^{\bullet}, \wedge, [\,,])$ over $\k$ 
is a $\Z$-graded supercommutative algebra $(\mathcal{V}, \wedge)$, together
with a bracket $[\,,]: \mathcal{V} \o \mathcal{V} \rightarrow \mathcal{V}$ of degree $-1$, such that the induced bracket of
degree zero on the shifted graded $\k$-module $\mathcal{V}^{\bullet+1}$ is
a Lie superbracket, satisfying the Leibniz identity,
\begin{equation}
[a \wedge b, c] = a \wedge [b, c] + (-1)^{m n} b \wedge [a, c], \quad a \in \mathcal{V}^m, b \in \mathcal{V}^n, c \in \mathcal{V}.
\end{equation}
\end{defn}

Finally, we recall the definition of a BV algebra.
\begin{defn}
  A Batalin-Vilkovisky (BV) algebra
  $(\mathcal{V}^\bullet, \wedge, \Delta)$ is a $\Z$-graded
  supercommutative algebra $(\mathcal{V}^\bullet, \wedge)$ equipped
  with an operator $\Delta: \mathcal{V} \rightarrow \mathcal{V}$ of
  degree $-1$ such that $\Delta^2 = 0$, and such that the bracket
  $[\,,]$ defined by
\begin{equation}\label{bvid}
(-1)^{m+1} [a,b] = \Delta(a \wedge b) - \Delta(a) \wedge b - (-1)^{m} a \wedge \Delta(b) + (-1)^{m+n} a \wedge b \wedge \Delta(1), \quad a \in \mathcal{V}^m, b \in \mathcal{V}^n,
\end{equation}
endows $(\mathcal{V}^\bullet, \wedge, [\,,])$ with a Gerstenhaber algebra
structure.  Here, $1 \in \mathcal{V}^0$ is the algebra unit.
\end{defn}

\section{General theory}\label{cyfrobs}
The goal of this section is to prove that periodic Calabi-Yau
Frobenius algebras have BV and Frobenius structures on their
Hochschild cohomology (Theorems \ref{cyfhhthm}, \ref{pcyfralg}).
Along the way, we will prove several general results about the
Hochschild (co)homology of Frobenius algebras over an arbitrary
base commutative ring.  We also explain why our
definition of Calabi-Yau Frobenius implies the CY condition of
\cite{ESk} (Theorem \ref{scyt}), and give a new proof that the
Hochschild cohomology of symmetric algebras is BV (\S
\ref{symfrobs}).

We will need a straightforward generalization of the stable module
category to a version relative to $\k$, which we relegated to Appendix
\ref{gfs}.

\subsection{Stable Hochschild (co)homology}
In this section, we define and begin the study of the \textbf{stable
  Hochschild (co)homology},
by replacing $\Ext$ by $\uExt$ in the definition. By Remark
\ref{taterem}, this is the Hochschild
version of Tate cohomology of finite groups.

For Hochschild homology, we will first need the notion of stable
tensor product:
\begin{defn}
Let $A$ be any algebra over $\k$ (projective as a symmetric $\k$-module).
For any $A^{\op}$-module $M$ and any $A$-module $N$,
such that $M,N$ are fg projective as $\k$-modules, define
\begin{multline}
  M \uo_A N := \{f \in M \o_A N \mid f \in \Ker(M \o_A N \rightarrow M
  \o_A I), \\ \text{ for any $\k$-split injection $N \into I$, with
    $I$ relatively injective}.\}
\end{multline}
\end{defn}
The condition that $f \in \Ker(M \o_A N \rightarrow M \o_A I)$ for any
particular $\k$-split injection $N \into I$ as above is equivalent to
the condition holding for \textit{all} such $\k$-split injections
(cf.~Appendix \ref{gfs}).

\begin{prop}
If $A$ is a Frobenius algebra, then
the definition of $\uo_A$ is symmetric in the following sense:
\begin{multline}
M \uo_A N \cong \{f \in M \o_A N \mid f \in \Ker(M \o_A N \rightarrow J \o_A N), \\ \text{ for any $\k$-split injection $M \into J$, with $J$ a 
relatively injective $A^\op$-module} \}.
\end{multline}
Thus, one has
\begin{equation}
M \uo_A N \cong N \uo_{A^\op} M.
\end{equation}
\end{prop}
\begin{proof}
Fix $\k$-split injections $N \into I, M \into J$, for 
$I, J$ relatively injective (=fg projective) $A$- and $A^\op$-modules,
respectively.  Since $I, J$ are projective, the maps
$M \o_A I \rightarrow J \o_A I$ and $J \o_A N \rightarrow J \o_A I$ are
injective.  Hence, the kernel of $M \o_A N \rightarrow J \o_A N$ is the
same as the kernel of $M \o_A N \rightarrow M \o_A I$ (both are the kernel
of $M \o_A N \rightarrow J \o_A I$).
\end{proof}
\begin{defn}
Let $A$ be a Frobenius algebra over $\k$. 
For any  $A^{\op}$-module $M$, and any $A$-module $N$, both
which are fg projective over $\k$, we define
the $i$-th stable Tor groups by
\begin{equation} \label{toridfn}
\uTor^A_i(M, N) :=  M \uo_A \Omega^i N.
\end{equation}
\end{defn}
\begin{prop} \label{torcor}
If $A$ is Frobenius over $\k$, then
\begin{equation}
\uTor^A_i(M, N) \cong \Tor^A_i(M, N), \quad i \geq 1,
\end{equation}
and moreover, the definition \eqref{toridfn} is symmetric in the following
sense:
\begin{equation} \label{toridfnsym}
M \uo_A \Omega^i N \cong \Omega^i M \uo_A N, \forall i \in \Z.
\end{equation}
\end{prop}
\begin{proof}
  The first statement follows similarly to Corollary \ref{extcor}.
  Similarly, we find that $\Omega^i M \uo_A N \cong \Tor_i(M, N)$ for
  $i \geq 1$, which yields \eqref{toridfnsym} for $i \geq 1$.  The
  statement is tautological for $i = 0$. To extend to negative $i$, we
  may use the trick
  $(\Omega^{-i} M) \uo_A \Omega^i (\Omega^{-i} N) \cong \Omega^i
  (\Omega^{-i} M) \uo_A \Omega^{-i} N$.
\end{proof}
\begin{defn} Suppose that $A$ is Frobenius over $\k$, and let $M$ be
  any $A$-bimodule which is fg projective as a
  $\k$-module.  The $i$-th stable Hochschild (co)homology groups
  (which only depend on the stable equivalence class of $M$) are
  defined by
\begin{equation}
\uHH^i(A, M) := \uExt^i_{A^e}(A, M), \quad \uHH_i(A, M) := \uTor_i^{A^e}(A, M).
\end{equation}
\end{defn}
\begin{cor} With $A, M$ as in the definition,
$\uHH^i(A, M) \cong \HH^i(A, M)$ and $\uHH_i(A, M) \cong \HH_i(A, M)$
for $i \geq 1$.
\end{cor}
\begin{rem}
  One could pose the definitions of stable Hochschild (co)homology and stable
  Ext and Tor when $A$ is not Frobenius, but one probably wants
  $\Omega$ to be an autoequivalence to have a reasonable notion (e.g.,
  if $A$ is ``relatively selfinjective;'' see Appendix \ref{gfs}).
\end{rem}

One has many algebraic structures attached to Hochschild cohomology
and homology: put together, these form the structure of calculus
(cf.~e.g., \cite{TT}; see Definition \ref{calcdfn} for the
definition). This includes cup products $\cup$ for Hochschild
cohomology, and contraction maps
$HH^j(A) \o HH_\ell(A) \rightarrow HH_{\ell-j}(A)$ for $j \leq \ell$.
For $f \in HH^j(A)$, we denote by
$i_f: HH_\ell(A) \rightarrow HH_{\ell-j}(A)$ the corresponding
contraction. We now show that the cup and contraction structures
extend to the stable, $\Z$-graded setting.
\begin{thm} \label{asfhh}
\begin{enumerate}
\item[(i)] Let $A$ be any Frobenius algebra over $\k$. Then one has a
  well-defined cup product on $\uHH^\bullet(A, A)$, giving the
  structure of an associative algebra, and extending the cup
product on $HH^{\geq 1}(A, A)$.
\item[(ii)]
One has a well-defined contraction operation $\uHH^{j}(A, A) \o \uHH_{k}(A,M) \rightarrow \uHH_{k-j}(A, M)$, which extends the usual contraction operation,
and satisfies the relation
\begin{equation} \label{contrint}
i_{f} i_{g} (x) = i_{f \cup g}(x),
\end{equation}
where $i_f(x)$ is the contraction of $f \in \uHH^\bullet(A, A)$ with
$x \in \uHH_\bullet(A, M)$. ($M$ is any fg $A$-bimodule.)
\item[(iii)] The algebra $\uHH^\bullet(A, A)$ is graded-commutative.
\end{enumerate}
\end{thm}
\begin{proof}
  (i) The cup product is easy to define: for $f \in \uHH^j(A, A)$ and
  $g \in \uHH^k(A, A)$, we have $f \in \uHom(\Omega^{j} A, A)$ and
  $\Omega^j g \in \uHom(\Omega^{j+k} A, \Omega^j A)$, so we may
  consider the composition
\begin{equation}
f \cup g := f \circ \Omega^j g = f \circ \Omega^{|f|} g,
\end{equation}
where $|f|$ denotes the Hochschild cohomology degree.
It follows immediately that the cup product is associative.  

(ii) To define the contraction operation, note that $f \in \uHH^j(A, A)
= \uHom(\Omega^j A, A)$ induces a map
\begin{equation}
\Omega^j A \uo_{A^e} \Omega^k A \rightarrow A \uo_{A^e} \Omega^k A,
\end{equation}
for all $k$.  Applying the equivalence $\Omega$, we obtain a map
\begin{equation}
A \uo_{A^e} \Omega^{j+k} A \rightarrow A \uo_{A^e} \Omega^k A,
\end{equation}
for all $k \in \Z$. This is the desired map.  We automatically get the
intertwining property \eqref{contrint}.  

(iii) We have isomorphisms in the stable module category,
\begin{equation}
\Omega^j A \o_A \Omega^k A \iso \Omega^{j+k} A,
\end{equation}
which follow from the fact that $\Omega^j A$ may be considered as
a projective left $A$-module for all $j$ (using the sequence
$\Omega^1 A \into A \o A \onto A$, which is split as left $A$-modules). We will need the 
\begin{claim} Let us use $\Omega^i A := (\Omega^1 A)^{\otimes_A i}$ for $i \geq 1$.
Given $f \in \uHom(\Omega^j A, \Omega^k A)$, we may form a representative of
$\Omega f \in \uHom(\Omega^{j+1} A, \Omega^{k+1} A)$ by either $\Id \o_A f$
or by $(-1)^{j-k} f \o_A \Id$.
\end{claim}
\begin{proof}
  This follows from the fact that the stable module category is a
  suspended category as in \cite{SA} (since it is a full monoidal
  subquotient of the derived category which is closed under
  suspension).  However, we give an explicit argument.  The terms of
  the normalized bar resolution may be written as
  $A \o \Omega^\bullet A$, and if we construct this by splicing
  together sequences
  $\Omega^{n+1} A \into A \o \Omega^{n} A \onto \Omega^n A$, it is
  easy to see that $f$ lifts to $f \o_A \Id$.  Writing the terms of
  the normalized bar resolution as $\Omega^\bullet A \o A$, we obtain
  the desired sign corrections of $(-1)^{j-k}$.
\end{proof}
As a consequence, if we have $f \in \uHH^j(A, A), g \in \uHH^k(A, A)$,
then we may compute $f \cup g$ in two ways. First, if $j, k \geq 0$,
then letting
$f' \in \uHom(\Omega^j A, A), g' \in \uHom(\Omega^k A, A)$, we may use
either the formula $f' \o_A g'$ or $(-1)^{|f| |g|} g' \o_A f'$.
Similarly, if $j, k$ are arbitrary, we take instead
$f' \in \uHom(\Omega^{a+j}A, \Omega^{a}A), g' \in \uHom(\Omega^{b+k}A,
\Omega^b A)$
for $a+j, a, b+k, b \geq 0$. Then the same argument for the
composition $\Omega^b f' \circ \Omega^{a+j} g'$ yields either
$(-1)^{j b} f' \o_A g'$ or
$(-1)^{(a+j)k} g' \o_A f'$. So, we obtain
$(-1)^{jb} f \cup g = (-1)^{(b+k)j} g \cup f$, as desired.
\end{proof}

Finally, we present a duality property for stable Hochschild homology,
which will induce a Frobenius algebra structure on Hochschild cohomology
for ``Calabi-Yau Frobenius algebras.'' To do this, we first need to
explain how to write standard complexes computing $\uHH_\bullet$ and
$\uHH^\bullet$. More generally, we define these computing stable Ext and Tor.
\begin{defn}
For any Frobenius algebra $A$ over $\k$ and any 
fg $A$-module $M$ which is projective over $\k$, call a \textbf{two-sided} resolution of $M$,
an exact $\k$-split complex of fg projective $A$-modules,
\begin{equation}
\cdots \rightarrow P_{2} \rightarrow P_{1} \rightarrow P_0 \rightarrow P_{-1} \rightarrow P_{-2} \rightarrow \cdots,
\end{equation}
such that $M$ is the cokernel of $P_1 \rightarrow P_0$ (and the
kernel of $P_{-1} \rightarrow P_{-2}$):
\begin{equation}
P_1 \rightarrow P_0 \onto M \into P_{-1} \rightarrow P_{-2}.
\end{equation}
\end{defn}
\begin{defn}
  For any Frobenius algebra $A$, any fg left $A$-module $M$ which is
  projective over $\k$, any fg right $A$-module $N$ which is
  projective over $\k$, and any two-sided resolution
  $P_\bullet$ of $M$, define the associated (``standard'') complex
  computing stable Tor,
\begin{equation}
C^{A}_\bullet(N, M) := N \o_A P_\bullet.
\end{equation}
Similarly, if now $M, N$ are both fg left $A$-modules which are
$\k$-projective, we define the associated (``standard'') complex
computing stable Ext,
\begin{equation}
C_A^\bullet(M, N) := \Hom_A(P_\bullet, N).
\end{equation}
We call the classes of the complexes $C^\bullet, C_\bullet$ in the
(unbounded) derived category of fg $\k$-modules, the
stable $M \uo_A^L N$ and $\uRHom(M, N)$.
\end{defn}
\begin{defn}
  For any two-sided $A^e$-resolution of $A$, we call
  $C_\bullet(A, A) := C^{A^e}_\bullet(A, A)$ and
  $C^\bullet(A, A) := C_{A^e}^\bullet(A, A)$ the associated
  ``standard'' complexes computing stable Hochschild homology and
  cohomology, respectively.
\end{defn}
\begin{rem}
  Note that we could have chosen to reserve the words ``standard'' for
  complexes resulting from the bar resolution of $A$ (which can be
  completed to a two-sided resolution, as we will explain).
\end{rem}
\begin{thm} \label{hhomduthm}
Let $A$ be a Frobenius algebra over $\k$. We have the following duality:
\begin{gather}
\Dk: \uHH_i(A, A) \iso \uHH_{-1-i}(A, A)^*, \quad \text{if $\k$ is a field},
\label{hhomdu} \\
\Dk: C^{A^e}_{\bullet}(A, A) \qiso C^{A^e}_{-1-\bullet}(A, A)^*, \quad \text{in general}.
\label{hchadu}
\end{gather}
Moreover, using \eqref{hhomdu}, the contraction maps become graded selfadjoint:
\begin{equation} \label{dkint}
\Dk(i_f x) = (-1)^{|f| |x|} i_f^* \Dk(x). 
\end{equation} 
\end{thm}
An easy extension of the theorem to coefficients in any bimodule $M$ which
is fg projective over $\k$ yields:
\begin{equation}
\D_\k: \uHH_i(A, M) \iso \uHH_{-1-i}(A, M^* \o_A A^\vee)^*, \quad
C^{A^e}_\bullet(A,A) \mathop{\iso}^{q.i.} C^{A^e}_{-1-\bullet}(A, M^* \o_A A^\vee)^*.
\end{equation}

To prove the theorem, the following easy identifications will be useful:
\begin{lemma} \label{eil}
Let $A$ be a Frobenius algebra over $\k$.
\begin{enumerate}
\item[(i)] For any left $A$-module $M$ which
is fg projective over $\k$,
\begin{equation}
  ({}_\phi M)^* \cong
  (M^*)_\phi, \quad ({}_\phi M)^\vee \cong 
  (M^\vee)_\phi, \quad {}_{\phi} A \o_A M \cong {}_\phi M.
\end{equation}
Furthermore,
\begin{equation} \label{mveeast}
M^\vee \o_A A^* \cong M^* \cong (M^\vee)_{\eta}.
\end{equation}
\item[(ii)] As $A$-bimodules, we have
\begin{equation} \label{avsenew}
_{\phi} A_\psi \cong A_{\phi^{-1} \circ \psi}, \quad A^\vee \cong
A_{\eta^{-1}}, \quad A^\vee \o_A A^* \cong A \cong A^* \o_A A^\vee.
\end{equation}
\end{enumerate}
\end{lemma}
\begin{proof}
(i) The first set of identities is immediate. For the second, we use
\eqref{sfe1}:
\begin{equation}
M^\vee \o_A A^*
\cong M^*_{\eta^{-1}} \o_A A^* \cong M^*.
\end{equation}
(ii) The first identity is clear. Applying
\eqref{mveeast} to the case of bimodules, we have
\begin{equation}
A^* \o_A A^\vee \o_A A^* \cong A^\vee \o_{A^e} (A^e)^* \cong A^*,
\end{equation}
which immediately gives the last identity, and hence the second.
\end{proof}
\begin{proof}[Proof of Theorem \ref{hhomduthm}]
Given any two-sided resolution $P_\bullet$ of $A$, the
lemma shows that the exact complex $P^\vee_{-1-\bullet} \o_A A^*$ must be
another two-sided resolution of $A$. Furthermore, 
we have
\begin{equation} \label{hhomisos}
(P^\vee_{-1-\bullet} \o_A A^*) \o_{A^e} A \cong P^\vee_{-1-\bullet} 
\o_{A^e} A^* \cong \Hom_{A^e}(P_{-1-\bullet}, A^*) \cong (A \o_{A^e} P_{-1-\bullet})^*,
\end{equation}
where the last isomorphism uses the standard adjunction.

To show the graded selfadjoint property \eqref{dkint}, we first note
the following naturality: applying $(\tilde f^\vee \o_A \Id_{A^*})
\o_{A^e} \Id$ on the LHS of \eqref{hhomisos} (where $\tilde f$ is the
lift of $f$ to $P_\bullet$) is the same as applying $(\Id \o_{A^e}
\tilde f)^*$ on the RHS.  Now, the LHS can be replaced by the complex
$\Hom_{A^e}(P_{-1-\bullet}, A^* \o_A P_{\bullet})$ (which is the total
complex, summing the $\bullet$ degrees).  Applying $(\tilde f^\vee
\o_A \Id_{A^*}) \o_{A^e} \Id$ becomes \textbf{right} composition with
$\tilde f$.  Using the same argument as in the proof of Theorem
\ref{asfhh}.(iii), this is chain-homotopic to applying \textbf{left}
composition with $(-1)^{|f| \cdot |-1-\bullet|} (\Id_{A^*} \o_A \tilde
f)$.  Now, via the quasi-isomorphism $A \o_{A^e} P_\bullet \simeq
(P^{\vee}_{-1-\bullet} \o_A A^*) \o_{A^e} P_{\bullet} \simeq
\Hom_{A^e}(P_{-1-\bullet}, A^* \o_A P_\bullet)$, applying $(-1)^{|f|
  \cdot |-1-\bullet|} \Id_{A^*} \o_A \tilde f$ gets carried to this
map, which we showed is chain-homotopic to the map obtained from
applying $(\Id \o_{A^e} \tilde f)^*$ on the RHS of \eqref{hhomisos}.
Now, passing to homology, we get the desired graded selfadjointness.
(In the Calabi-Yau Frobenius case which we will define, we can deduce
this result more simply as in the proof of Theorem \ref{cyfhhthm}).
\end{proof}

\subsection{Relative Serre duality for the stable module category}
From now on, let $A$ be a Frobenius algebra over $\k$.  As was noticed
in \cite{ESk}, when $\k=$ a field, the Auslander-Reiten homomorphisms
give a Serre duality for the stable module category. We recall this
and produce a relative version.
\begin{ntn}
Let $\nu: A-mod \rightarrow A-mod$ be the \textbf{Nakayama functor}, $\nu = * \circ \vee$.
\end{ntn}
Note that, when $\k$ is a field, $\nu$ sends projectives to injectives,
and vice-versa, and in fact induces an equivalence of categories
$\{\text{projective $A$-modules}\} \leftrightarrow \{\text{injective
  $A$-modules}\}$.   Also,
by \eqref{sfe1}, one has the following simple formula for $\nu M$ (which
will be useful later):
\begin{equation} \label{nuform}
\nu M \cong A^* \o_A M \cong {}_{\eta^{-1}} M.
\end{equation}

\begin{prop}\cite{ESk}
When $\k$ is a field, we have functorial isomorphisms
\begin{equation}
\uHom_A(M,N) \cong \uHom_A(N, \Omega \nu M)^*.
\end{equation}
\end{prop}
The proof is based on the Auslander-Reiten formulas (cf.,
e.g.~\cite{ARS,ASS}).  

In order to make proper sense of the Serre duality for arbitrary $\k$,
it is necessary to replace the groups above by complexes and dual
complexes.  First, observe that the above isomorphism
can be rewritten, by replacing $M$ by $\Omega^i M$, as
\begin{equation}
\uExt^i_A(M, N) \cong \uExt^{-1-i}_A(N, \nu M)^*.
\end{equation}
The following then gives a version of the above for general $\k$:
\begin{thm}\label{comt} 
Let $M, N$ be
$A$-modules which are fg projective as $\k$-modules.
Then, one has a functorial quasi-isomorphism in the derived category,
\begin{equation} \label{comte}
C_A^\bullet(M, \nu N) \simeq C_A^{-1-\bullet}(N, M)^*.
\end{equation}
\end{thm}
\begin{proof}
  Fix a two-sided projective resolution $P_\bullet$ of $N$.  Applying
  $\nu$, we obtain a resolution of $\nu(N)$.  By Proposition
  \ref{riap}, $\Hom_A(M, \nu P_{-1-\bullet})$ may be used to compute
  $\Ext_A^\bullet(M, \nu N)$, since $\nu P_\bullet$ is a two-sided
  resolution of $\nu N$ consisting of relatively injectives.
  Furthermore, for any two-sided resolution $Q_\bullet$ of $M$, we may
  obtain quasi-isomorphisms
\begin{equation}
  C_A^\bullet(M, \nu N) = \Hom_A(Q_{\bullet}, \nu N) \liso 
  \Hom_A(Q_{\bullet}, \nu P_{-1-\bullet}) \iso \Hom_A(M, \nu P_{-1-\bullet}).
\end{equation}
Next,  we show that
  $\Hom_A(M, \nu P_{-1-\bullet}) \simeq \Hom_A(P_{-1-\bullet}, M)^*$.  For any module $L$, there
  is a functorial map $L^\vee \o_A M \rightarrow \Hom_A(L, M)$, which is an isomorphism if $L$ is projective.
  Applying this map and its dual, we obtain
\begin{equation}
\Hom_A(P_{-1-\bullet}, M)^* \iso (P_{-1-\bullet}^\vee \o_A M)^* \mathop{\cong}^{adj.} \Hom_A(M, \nu P_{-1-\bullet}),
\end{equation}
the last map using adjunction.
\end{proof}

\subsection{Calabi-Yau and periodic Frobenius algebras}
We will modify the definition of Calabi-Yau algebra \cite{Gcy} to
suit Frobenius algebras.
First, we recall this definition.
\begin{defn}\cite{Gcy} An associative algebra $A$ over a commutative
  ring $\k$ is called
  \textbf{Calabi-Yau} if it has finite Hochschild dimension, and one
  has a quasi-isomorphism in the derived category of $A^e$-modules,
\begin{equation}
f: A[d] \iso \RHom_{A^e}(A, A \o A)
\end{equation}
which is self-dual:
\begin{equation} \label{selfdual}
f^! \circ \iota = f[-d],
\end{equation}
where, for any map $g: M \rightarrow N$ of $A$-bimodules,
$g^!: \RHom_{A^e}(N, A^e) \rightarrow \RHom_{A^e}(M, A^e)$ is the
natural map, and
$\iota: A \rightarrow \RHom_{A^e}(A, \RHom_{A^e}(A, A \o A))$ is the
natural map. Here $[-]$ denotes the shift in the derived category.
\end{defn}

However, if $A$ is Frobenius and of finite Hochschild dimension (part of being
(usual) CY), then $\Omega^i A = 0$ for large enough $i$, and hence $A \simeq 0$ in the stable module category. Hence, $HH^i(A, M) = 0$ for all $i \geq 1$, 
so $A$ has
Hochschild dimension zero. That is, $A$ is a projective $A$-bimodule (i.e.,
$A$ is separable).

For a separable Frobenius algebra to additionally be Calabi-Yau of
dimension zero, we require exactly that $A \cong A^\vee$ as
$A$-bimodules. By Lemma \ref{eil}, this is equivalent to $A^* \cong
A$, i.e., $A$ must be a symmetric, separable Frobenius algebra.  

This is not general enough, so we replace this notion.
Let us first restate
the usual Calabi-Yau property in terms of a quasi-isomorphism of complexes:

\begin{equation}
\xymatrix{
0 \ar[r] & P_m \ar[r] \ar[d] & P_{m-1} \ar[r] \ar[d] & \cdots \ar[r] \ar[d] & P_0 \ar@{->>}[r] \ar[d] & A \ar@{=}[d] \\
0 \ar[r] & P_0^{\vee} \ar[r] & P_1^{\vee} \ar[r] & \cdots \ar[r] & P_m^{\vee} \ar@{->>}[r] & A.
}
\end{equation}
Here and from now on, the functor $^\vee$ will be in the category of
$A$-bimodules, i.e., $M^\vee = \Hom_{A^e}(M, A \o A)$.

If now $A$ becomes a Frobenius algebra, as explained earlier, we
cannot have such a quasi-isomorphism (for $m \geq 1$): in fact, when
dualizing the top sequence, we get something that begins with
$A^\vee \into P_0^\vee \rightarrow P_1^\vee$.  However, 
it still makes sense to ask for
a commutative diagram, with exact rows, as follows:
\begin{equation}
  \xymatrix{
    A^{\vee} \ar@{^{ (}->}[r] & P_{m} \ar[r] \ar[d] & P_{m-1} \ar[r] 
\ar[d] & \cdots \ar[r]\ar[d] & P_0 \ar[d] \ar@{->>}[r] & A \ar@{=}[d] \\
    A^\vee \ar@{^{ (}->}[r] & P_0^{\vee} \ar[r] & P_1^{\vee} \ar[r] & 
\cdots \ar[r] & P_m^{\vee} \ar@{->>}[r] & A
  }
\end{equation}
In fact, since the dual is automatically exact in the Frobenius case, by
a standard result of homological algebra, such a diagram must
automatically exist given a resolution as in the top row. So it is
enough to ask for such a resolution.

On the level of the stable module category, having such a resolution
implies the following condition:
\begin{defn}
A Frobenius algebra $A$ over $\k$ is called \textbf{Calabi-Yau Frobenius
of dimension $m$} if
one has isomorphisms in the stable module category $\Stab_\k(A^e)$, for
some
$m \in \Z$:\footnote{The self-duality property 
required for Calabi-Yau algebras \eqref{selfdual} can still be imposed
here: this would be the condition that
$f: A^\vee \simeq \Omega^{m+1} A$ satisfies
$f^\vee = \Omega^{-m-1} f$. 
We do not need this.}
\begin{equation} \label{cyfrob}
f: A^\vee \simeq \Omega^{m+1} A.
\end{equation}
If there is more than one such $m$, then we pick the smallest nonnegative value of $m$
(which exists because such algebras must be periodic as in the subsequent definition).

If additionally $A$ has a grading, such that the above isomorphism
is a graded isomorphism when composed with some shift (considering the
stable module category to now be graded), then we say that $A$ is a
\textbf{graded Calabi-Yau Frobenius algebra}. More precisely, if
$f: A^\vee(m') \simeq \Omega^{m+1} A$ 
is a graded isomorphism, where $(\ell)$ denotes
the shift by $\ell$ with the new grading, then one says that $A$ is \textbf{graded Calabi-Yau 
Frobenius with dimension $m$ of shift $m'$.}
\end{defn}
We remark that the above definition of Calabi-Yau Frobenius is
(apparently) stronger than the notion of Calabi-Yau for selfinjective
algebras discussed in \cite{ESk} for the case $\k$ is a field: see
Theorem \ref{scyt} and the comments thereafter.

\begin{defn}
A Frobenius algebra over $\k$ is called \textbf{periodic Frobenius of
period $n$} for some $n > 0$ if one has (in $\Stab_\k(A^e)$)
\begin{equation}\label{perfrob}
g: A \simeq \Omega^n A,
\end{equation}
and that is the smallest positive $n$ for which one has such an isomorphism.
If $A$ has a grading, we define \textbf{graded periodic Frobenius of period $n$ and
shift $n'$} as before (if $g: A(n') \simeq \Omega^n A$ is a graded isomorphism in
the stable module category).
\end{defn}
Note that it makes sense to be Calabi-Yau Frobenius of
\textbf{negative} dimension.  In particular, \textit{any symmetric
  Frobenius algebra is either Calabi-Yau Frobenius of dimension $-1$,
  or periodic Calabi-Yau Frobenius of dimension $n-1$ and period $n$
  for some $n \geq 1$.}

Also, any periodic Frobenius algebra must have even period, unless
$2 \cdot \Id \simeq 0$ in the stable module category, e.g.,
$\text{char } \k = 2$, as we will see in Theorem \ref{perfhhthm}.  
In
particular, the CY dimension must be odd for symmetric Frobenius algebras.

\begin{exam} \label{dynqe} The preprojective algebras of ADE Dynkin
  quivers are periodic Calabi-Yau Frobenius (using \cite{EE2}) of
  dimension $2$ and shift $2$, and of period $6$ and shift $2h$ (twice
  the Coxeter number).  The essential ingredient is the Schofield
  resolution \cite{RS} (cf., e.g., \cite{EE2}), with $R = \k^I$, where $I$ is the vertex set:
\begin{equation}
A^\vee(2) \into A \o_R A(2) \rightarrow A \o_R V \o_R A \rightarrow A \o_R A \onto A.
\end{equation}
Here, $V$ is the free $\k$-module spanned by the edges of the quiver. This
completes to a periodic projective resolution of length $6$, since the
Nakayama automorphism has order $2$ (for more details, see Proposition
\ref{fonpfr}).  The fact that the shift is $2h$ follows from the fact
that $A^\vee(2) \cong A_{\eta}(h)$, which amounts to the fact that the
degree of the image of $\Id \in A \o_R A^*$ under the isomorphism
$A \o_R A^* \cong A \o_R A$ from the pairing, is $h-2$ (i.e., $h-2$ is the
degree of the product of any basis element with its dual basis
element).

We note that an important part of the above is showing that $\Pi_Q$ is
free (equivalently, projective) over $\Z$; this follows from explicit
$\Z$-bases (such as those in \cite{Eu2}, which one may verify is
integral; types $A,D$ over $\Z$ are also in \cite[\S 4.2.2]{S3}).  
\end{exam}
\begin{exam} \label{cedynqe}
Similarly, the centrally extended preprojective algebras \cite{ER},
over $\k=\C$, are periodic Calabi-Yau Frobenius with dimension $3$ (of
shift $4$) and period $4$ (of shift $2h$), using \cite{Eu}.  These
are symmetric, i.e., have trivial Nakayama automorphism.  Note that
these are not, in general, torsion-free over $\Z$ (cf.~\S
\ref{ces}), hence not Frobenius over $\Z$, although the definition
may be modified to correct this.
\end{exam}
\begin{exam}\label{gdynqe}  Similarly to the Dynkin case, one may
  consider preprojective algebras of \textit{generalized} Dynkin type:
  this refers to preprojective algebras of type $T_n$ (otherwise known
  as $L_n$) which can be obtained from $\Pi_{A_{2n}}$ by passing to
  fixed points under the Nakayama automorphism,
  $\Pi_{T_n} := (\Pi_{A_{2n}})^{\eta}$.  In other words, this is associated
  to a graph of $T_n$ type.

In \cite{Eu4} (using a variant of the Schofield resolution \cite{RS} (cf.~\cite{BES}),
it is proved that $\Pi_{T_n}$ is periodic Calabi-Yau Frobenius of
dimension $5$ (and shift $h+2$) and period $6$ (and shift $2h$).  Also,
in \cite{Eu4}, a correction to results of \cite{ES2} in type $A$ is given.
\end{exam}
\begin{exam}\cite{BBK}
  The trivial extension algebras of path algebras of Dynkin quivers
  are periodic Calabi-Yau Frobenius (in fact, symmetric) of dimension
  $2h-3$ and period $2h-2$, where $h$ is the Coxeter number.  These are ``almost-Koszul dual'' to the preprojective algebras; see \cite{BBK}.
\end{exam}
For Calabi-Yau Frobenius algebras, being periodic Frobenius
is closely related to having finite-order Nakayama automorphism:
\begin{prop} \label{fonpfr} If $A$ is Calabi-Yau Frobenius of dimension
$\neq -1$, then the following
are equivalent:
\begin{enumerate}
\item[(i)] $A$ is periodic Frobenius;
\item[(ii)] For some $p > 0$, one has
\begin{equation} \label{fonak}
A \simeq A_{\eta^p} \text{ in the stable bimodule category}\footnote{Note that, if we had defined Calabi-Yau and periodic Frobenius algebras using complexes rather than the stable module category, then this would be replaced by an honest isomorphism of bimodules,
$A \cong A_{\eta^p}$, and
hence $\eta^p$ would  have to be inner.}
\end{equation}
(e.g., if $\eta^p$ is inner).
\end{enumerate}
In the situation that the above are satisfied, then the Calabi-Yau
dimension $m$, the period $n$, and the smallest $p > 0$ such that
\eqref{fonak} holds, are related by
\begin{equation} \label{mpn}
n = p \cdot \gcd(n, m+1),
\end{equation}
and $r = \gcd(n, m+1)$ is the smallest positive integer such that $\Omega^r A \simeq A_{\eta^k}$ for some integer $k$.
\end{prop}
In particular, a Calabi-Yau Frobenius algebra of dimension $\leq -2$
must have infinite-order Nakayama automorphism (since if it were periodic,
the CY dimension is defined to be nonnegative). Being dimension $-1$ is a special
case, consisting of stably symmetric Frobenius algebras that are not periodic.
\begin{proof} (ii) implies (i):  We have $\Omega^{m+1} A \simeq A^\vee \cong A_{\eta^{-1}}$.  For any $p$ such that \eqref{fonak} holds, we have $A \simeq A_{\eta^{-p}} \simeq \Omega^{p \cdot (m+1)} A$, yielding (i). 

(i) implies (ii):  As in the previous paragraph, we deduce that $\Omega^{k \cdot (m+1)} A \simeq A_{\eta^{-k}}$.  If $\Omega^n A \simeq A$, then we would
deduce that $A \simeq \Omega^{-n \cdot (m+1)} A \simeq A_{\eta^{n}}$, yielding (ii).

To obtain \eqref{mpn}, let $r > 0$ be the smallest positive integer
such that $\Omega^r A \simeq A_{\eta^k}$ for some integer $k$.  It
follows that $r \mid (m+1)$ and $r \mid n$.  Since $r \mid (m+1)$, it
must be that $k$ is relatively prime to $p$, otherwise
$A_{\eta^{k \cdot \frac{m+1}{r}}} \simeq \Omega^{m+1} A \simeq
A_{\eta^{-1}}$
would contradict the minimality of $p$.  Similarly, we deduce that
$n = p \cdot r$, and $r = \gcd(n,m+1)$.
\end{proof}
Also, we have the following growth criterion:
\begin{ntn} For any fg $\k$-module $M$, let
  $g(M) = g_\k(M) \geq 0$ denote its minimal number of generators.
\end{ntn}
\begin{prop}
  Suppose that $\k$ has finite Krull dimension (or just that its
  maximal ideal spectrum has finite dimension).  If $A$ is Calabi-Yau
  Frobenius of dimension $\neq -1$, or if $A$ is periodic Frobenius,
  then for any fg $A$-modules $M, N$ which are
  projective as $\k$-modules, there is a positive integer $p \geq 1$
  such that, for all $i \in \Z$, $\uExt^i(M, N)$ and $\uTor^i(M^*, N)$
  are generated by at most $p$ generators over $\k$. We may take
  $p = a g(M) \cdot g(N) + b$ for some $a, b \geq 0$ depending
  only on $A$ and $\k$.

  Moreover, for any fg $A$-bimodule $L$ which is
  projective as a $\k$-module, there is a positive integer $q \geq 1$
  such that $\uHH^i(A, L), \uHH_i(A, L)$ are generated by at most $q$
  generators, for all $i \in \Z$. Again, we may take $q = a' \cdot g(L) + b'$
for some $a', b' \geq 0$ depending only on $A$ and $\k$.
\end{prop}
\begin{proof} Without loss of generality, assume that the 
tautological map $\k \rightarrow A$ is injective. Under either assumption
 of the proposition ($A$ is CY Frobenius of dim $\neq -1$ or $A$ is periodic
Frobenius), we have $\Omega^r A \simeq A_\phi$ for some
  $r \geq 1$ and some automorphism $\phi$ of $A$.  Thus, for all
  integers $s$, we may write $s = r \cdot k + r'$ for some
  $0 \leq r' < r$, and then
  $\Omega_{A^e}^s A \simeq \Omega_{A^e}^{r'} A_{\phi^k}$.  
Since $\Omega_A^s M
  \simeq \Omega_{A^e}^s A \o_A M$, we have
\begin{equation}
\uExt^s(M, N) \cong \uHom(\Omega_{A^e}^s A \o_A M, N) \cong \uHom({}_{\phi^{-k}} \Omega_{A^e}^{r'}A \o_A M, N).
\end{equation}
Furthermore, using the normalized bar resolution of $A$ (twisted by automorphisms of $A$), we may take $\Omega_{A^e}^s A$ to be isomorphic, as a right $A$-module, to $(A/\k)^{\o s} \o A$.  Then, localizing $\k$ at any prime ideal,
$\uHom({}_{\phi^{-k}} \Omega_{A^e}^{r'}A \o_A M, N)$ can have rank at most
$(g(A)-1)^{r'} \cdot g(M) \cdot g(N)$.  By Theorem 1 of \cite{Sw2},
$\uHom({}_{\phi^{-k}} \Omega_{A^e}^{r'}A \o_A M, N)$ can then be
generated by at most $(g(A)-1)^{r'} \cdot g(M) \cdot g(N) + d$
elements, where $d$ is the dimension of the maximal ideal spectrum of
$\k$.

We can apply the same idea for $\uTor^s(M, N)$.  For $\uHH^j(A, L)$ and
$\uHH_j(A, L)$ where $L$ is now an $A$-bimodule, we now use that these
are given by $\uHom_{A^e}(\Omega^j A, L)$ and $\Omega^j A \uo_{A^e} L$,
and apply the same argument as above.
\end{proof}
In particular, the above proposition rules out tensor products of
periodic Frobenius algebras over a field (which are not separable)
from being periodic or Calabi-Yau Frobenius of dimension $\neq -1$, by
the K\"unneth theorem.

The following theorems also have graded versions (by incorporating
the shifts in the definitions), but we omit them for simplicity.
\begin{thm} \label{scyt}
\begin{enumerate}
\item[(i)] For any Calabi-Yau Frobenius algebra $A$ of dimension $m$,
  and any fg $A$-modules $M, N$ which are projective over $\k$, one
  has isomorphisms (functorial in $M, N$)
\begin{gather}
\uHom_A(M, N) \cong \uHom_A(N, \Omega^{-m} M)^*, \quad \text{if $\k$ is a field}, \\
\uRHom^\bullet_A(M, N) \simeq \uRHom^{m-\bullet}_A (N, M)^* \quad \text{generally.}
\end{gather}

That is, $\Omega^{-m}$ is a (right) Serre
functor for the stable module category relative to $\k$.
\item[(ii)] For any periodic Frobenius algebra $A$ of period $n$, and any
fg $A$-modules $M, N$, one has isomorphisms
\begin{equation}
\uHom_A(M, N) \cong \uHom_A(\Omega^n M, N), \quad  \uRHom^\bullet_A(M, N) \simeq
\uRHom^{\bullet+n}_A(M, N).
\end{equation}
\end{enumerate}
\end{thm}
As a corollary, we see that any Calabi-Yau Frobenius algebra over a
field is also Calabi-Yau in the sense of \cite{ESk}. We do
not know if the reverse implication holds.
\begin{proof}
(i) By \eqref{nuform} (which holds for arbitrary $\k$) and \eqref{avsenew},
we have a stable equivalence
\begin{equation}
M \simeq \Omega^{m+1} \nu(M), \text{ i.e., } \nu M \simeq \Omega^{-m-1} M.
\end{equation}
Since $\Omega \nu$ already provides a right Serre functor in the
senses needed for (i) (using Theorem \ref{comt}), we now know that
$\Omega^{-m}$ also provides such a right Serre functor (note that the
above stable equivalences are clearly functorial). 

(ii) This is easy.
\end{proof}
\begin{rem}\label{twistingrem}
  For general Frobenius algebras, one can still say that
  $\uHom(M, N) \cong \uHom(N, A^* \o_A \Omega M)^*$ (by \eqref{nuform}
  and Theorem \ref{comt}), so that the Serre functor involves shifting
  and twisting by an invertible (under $\o_A$) bimodule. Similarly,
  the following results on Hochschild (co)homology have analogues for
  arbitrary Frobenius algebras using twisting (by powers of $\eta$) as
  well as shifting. For this, one considers the bigraded algebra
  $\uHom^{\bullet}_{A^e}(A, A_{\eta^{\bullet'}})$ and its bigraded
  module $\Omega^{\bullet} A \uo_{A^e} A_{\eta^{\bullet'}}$.  We do
  not need this for our examples.  However, it might be interesting to
  try to apply this formalism to finite-dimensional Hopf algebras
  (analogously to \cite{BZ}).
\end{rem}
We now present results on stable Hochschild cohomology of
periodic and CY Frobenius algebras:
\begin{thm} \label{cyfhhthm}
\begin{enumerate}
Let $A$ be a Calabi-Yau Frobenius algebra of dimension $m$, and $M$
any $A$-bimodule.
\item[(i)] 
One has isomorphisms and quasi-isomorphisms
\begin{equation}
\D: \uHH^\bullet(A, M) \cong \uHH_{m-\bullet}(A, M), \quad
\uRHom^\bullet_{A^e}(A, M) \simeq (A \uo^L_{A^e} M)_{m-\bullet}. \label{hhehh}
\end{equation}
For $M=A$, these isomorphisms intertwine cup product with contraction:
\begin{equation}\label{hhhhint}
\D \circ i_f(x) = f \cup \D(x), \quad \forall x \in \uHH_\bullet(A,A), f \in \uHH^\bullet(A, A).
\end{equation}
\item[(ii)] Let $M$ be fg projective over $\k$. One has
\begin{gather} \label{hhdual}
\uHH^\bullet(A,M) \cong \uHH^{m-\bullet}(A, M^*)^*, \quad \text{if $\k$ is
a field,} \\ \label{rhhdual}
\uRHom^\bullet_{A^e}(A, M) \simeq \uRHom^{m-\bullet}_{A^e}(A, M^*)^*, \quad \text{generally.}
\end{gather}
In the case $M = A$, we may rewrite this, respectively, as
\begin{equation} \label{hhduals}
\uHH^\bullet(A,A) \cong \uHH^{(2m+1)-\bullet}(A, A)^*, \quad
\uRHom^\bullet_{A^e}(A, A) \simeq \uRHom^{(2m+1)-\bullet}_{A^e}(A, A)^*.
\end{equation}
\item[(iii)] The induced pairing
\begin{equation}
\uHH^{\bullet}(A, A) \o \uHH^{(2m+1)-\bullet}(A, A) \rightarrow \k
\end{equation}
is invariant with respect to cup product:
\begin{equation} \label{cupinv}
(f, g \cup h) = (f \cup g, h), \quad |f|+|g|+|h| = 2m+1,
\end{equation}
and is nondegenerate if $\k$ is a field.
Moreover, for all $\k$, one has a nondegenerate invariant pairing
in the derived category (of degree $-(2m+1)$),\footnote{Here, a pairing $X \o Y \rightarrow \k$ is nondegenerate if it induces a (quasi-)isomorphism $X \simeq Y^*$, or equivalently, $Y \simeq X^*$.}
\begin{equation}
\uRHom_{A^e}(A, A) \o \uRHom_{A^e}(A, A) \rightarrow \k,
\end{equation}
inducing \eqref{hhduals}.

In other words, $C^\bullet$ is a Frobenius algebra in the derived
category of $\k$-modules, and if $\k$ is a field, $\uHH^\bullet$ is a
graded Frobenius algebra over $\k$ (using a definition
that only requires finite-generation in each degree).  
\end{enumerate}
\end{thm}
We also remark that \eqref{cupinv} and \eqref{hhhhint}, together with the
graded commutativity of Theorem \ref{asfhh}.(iii), give another proof
of the graded selfadjointness of $i_-^*$ \eqref{dkint} in this case.
\begin{proof}
(i) Let us pick a two-sided resolution of $A$:
\begin{equation} \label{tsra}
\cdots \rightarrow P_{2} \rightarrow P_1 \rightarrow P_0 \onto A \into P_{-1} \rightarrow P_{-2} \rightarrow \cdots,
\end{equation}
so that, removing the $A$, the complex $P_\bullet$ is an exact complex
of projectives. Since $\Omega^{m+1} A \simeq A^\vee$ (in the stable
module category), we may form chain maps between the two obtained
resolutions,
\begin{equation}
P_{\bullet+(m+1)} \leftrightarrow P_{-1-\bullet}^\vee,
\end{equation}
such that their composition on either side induces identity maps on the level
of $\uHom(\Omega^i A, \Omega^i A)$ and $\uHom(\Omega^i A^\vee, \Omega^i A^\vee)$ for all $i \in \Z$.  As a result, upon applying the functor $A \o_{A^e} -$, 
we obtain quasi-isomorphisms of the resulting complexes,
\begin{equation} \label{cydpf}
P_{\bullet+(m+1)} \o_{A^e} M \iso P_{-1-\bullet}^\vee \o_{A^e} M.
\end{equation}
However, it is clear that these are standard complexes computing $\uHH_{m+1+\bullet}(A, M)$ and
$\uHH^{-1-\bullet}(A, M)$, which is what we needed.

To prove the intertwining property \eqref{hhhhint}, we note that, for
$f \in \uHH^j(A, A)$, applying $i_f$ in the LHS of \eqref{hhhhint} is
the same as applying the corresponding element of
$\uHom(A, \Omega^{-j} A)$ to the $M=A$ in \eqref{cydpf}.  
Similarly, applying
$f \cup -$ to the RHS is post-composing with $f$, which is
applying the same element of $\uHom(A, \Omega^{-j} A)$ to the $A$ in the
RHS of \eqref{cydpf}.

(ii) This follows from part (i) and Theorem \ref{hhomduthm}. To fix
the signs (so as to obtain the Frobenius property in part (iii)), we use for
the duality
$x \mapsto (-1)^{|x| \cdot m} \D^* \circ \D_\k \circ \D(x)$, rather
than only $\D^* \circ \D_\k \circ \D$.
Alternatively, we can use
$\nu A = A^{\vee *} \simeq A_{\eta^2} \simeq \Omega^{-2(m+1)}A$ and
Theorem \ref{comt}, and similarly for arbitrary $M$, which will give
the desired property.

(iii) For this, we will use functoriality of the isomorphisms in (ii). 
(The result also follows from Theorem \ref{hhomduthm}.(iii).)

In the case that $\k$ is a field, the first isomorphism of
\eqref{hhduals} comes from the functorial isomorphisms of
\eqref{comte} (in the case $\k$ is a field, we can rewrite this as
$\uHom(M, N) \iso \uHom(N, \nu \Omega M)^*$).  Thus, we have the
commutative square (for any $g \in \uHH^{j}(A, A)$)
\begin{equation}
\xymatrix{
\uHom_{A^e} (A, \Omega^{-i} A) \ar[d]^{\Omega^{-i-j} g\, \circ} 
\ar[r]^-\sim & \uHom_{A^e}(\Omega^{-i} A, \Omega^{-(2m+1)} A)^* \ar[d]^{(\circ\,\Omega^{-i} g)^*}  \\
\uHom_{A^e} (A, \Omega^{-i-j} A) \ar[r]^-\sim & 
\uHom_{A^e}(\Omega^{-i-j} A, \Omega^{-(2m+1)} A)^*
}
\end{equation}
which gives exactly the invariance property needed.  On the level of
complexes, we may make the desired statement as follows.
Let $P_\bullet$,
as in \eqref{tsra}, be a two-sided resolution of $A$. Then, we have the
following sequence, where the isomorphisms mean quasi-isomorphisms,
and all complexes are total complexes graded by $\bullet$:
\begin{equation}
C^\bullet \liso \Hom_{A^e}(P_{\bullet}, P_{-\bullet}) \iso 
\Hom_{A^e}(P_{-\bullet}, \nu P_{1+\bullet})^* \iso \Hom_{A^e}(P_{-\bullet}, P_{\bullet-(2m+1)})^* \iso (C^{(2m+1)-\bullet})^*.
\end{equation}
Next, for any $j \in \Z$ and any $g \in \Hom_{A^e}(P_{\bullet+j}, P_\bullet)$,
we have the commutative square
\begin{equation}
\xymatrix{
\Hom_{A^e} (P_{\bullet}, P_{-\bullet}) \ar[d]^{g\,\circ} 
\ar[r]^-\sim & \Hom_{A^e}(P_{-\bullet}, P_{\bullet-(2m+1)} A)^* \ar[d]^{(\circ\, g)^*}  \\
\Hom_{A^e} (P_{\bullet}, P_{-\bullet-j}) \ar[r]^-\sim & 
\Hom_{A^e}(P_{-\bullet-j}, P_{\bullet-(2m+1)})^*
}
\end{equation}
which gives the desired result. The pairing
$C^\bullet \o C^\bullet \rightarrow \k$ in the derived category is
given by replacing $C^\bullet$ with
$\Hom_{A^e} (P_{\bullet}, P_{-\bullet})$ and using
composition.
\end{proof}
We now need to recall the definition of calculus:
\begin{defn}\cite{CST} (p. 93) \label{calcdfn}
A precalculus is a pair of a Gerstenhaber algebra $(\mathcal{V}^\bullet, \wedge, [\,,])$ and a graded vector space $\mathcal{W}^{\bullet}$  together with
\begin{enumerate}
\item A module structure $\iota_{-}: \mathcal{V}^\bullet \o \mathcal{W}^{-\bullet} \rightarrow \mathcal{W}^{-\bullet}$ of the graded commutative algebra $\mathcal{V}^\bullet$ on $\mathcal{W}^{-\bullet}$;
\item An action $\mathcal{L}_{-}$ of the graded Lie algebra $\mathcal{V}^{\bullet+1}$ on $\mathcal{W}^{-\bullet}$, which satisfies the following compatibility
conditions:
\begin{gather} \label{compat1}
\iota_a \mathcal{L}_b - (-1)^{|a|(|b|+1)} \mathcal{L}_b \iota_a = \iota_{[a,b]}, \\
\label{compat2} \mathcal{L}_{a \wedge b} = \mathcal{L}_a \iota_b + (-1)^{|a|} \iota_a \mathcal{L}_b.
\end{gather}
\end{enumerate}
A calculus is a precalculus
$(\mathcal{V}^\bullet, \mathcal{W}^\bullet, [\,,], \wedge, \iota_-,
\mathcal{L}_-)$
together with a differential $d$ of degree $1$ on
$\mathcal{W}^\bullet$ satisfying the Cartan identity:
\begin{equation} \label{cartan}
\mathcal{L}_a = d \iota_a - (-1)^{|a|} \iota_a d.
\end{equation}
\end{defn}
It is a result of \cite{DGT} that, for any associative algebra $A$, the collection
\begin{equation}
  (\HH^\bullet(A,A), \HH_\bullet(A, A), \{,\}, \cup, i_-, \mathcal{L}_-, B)
\end{equation} 
is a calculus, where $B$ is the Connes differential, $\{,\}$ the
Gerstenhaber bracket, and $\mathcal{L}_-$ the Lie derivative
operation.
\begin{thm} \label{perfhhthm}
Let $A$ be any periodic Frobenius algebra of period $n$.  Then, $n$
must be even if $2 \cdot \uHH^\bullet(A,A) \neq 0$, and
\begin{enumerate}
\item[(i)] One has isomorphisms and quasi-isomorphisms
\begin{gather}
  \uHH^\bullet(A, A) \cong \uHH^{\bullet+n}(A, A), \quad
  \uHH_\bullet(A, A) \cong \uHH_{\bullet+n}(A, A), \label{uhhper} \\
  C_{A^e}^\bullet(A,A) \iso C_{A^e}^{\bullet+n}(A,A),
\end{gather}
and similarly a quasi-isomorphism between the standard complex
computing stable Hochschild homology and its shift by $n$.  Moreover,
the isomorphisms may be induced by cup product on the left 
with the element $1' \in \uHH^n(A, A)$ representing the
given isomorphism $\Omega^n A \iso A$, and by the contraction $i_{1'}$.
\item[(ii)] The stable Hochschild cohomology is a Gerstenhaber algebra, which
extends the Gerstenhaber structure on usual Hochschild cohomology;
\item[(iii)] The stable Hochschild cohomology and stable Hochschild homology form the structure of a calculus, extending the usual calculus structure.
\end{enumerate}
\end{thm}
\begin{proof}
(i) The isomorphisms follow as in previous proofs from the stable module isomorphism
$1': A \simeq \Omega^n A$.  To show that they are induced by cup product or contraction with $1' \in \uHH^n(A, A)$, let us construct a projective resolution\footnote{This is isomorphic to the normalized bar resolution.} of $A$
such that 
\begin{equation}
\Omega^j A \cong (\Omega^1 A)^{\o_A j}.
\end{equation}
Following the proof of Theorem \ref{hhomduthm}.(iii), we construct
this from any exact sequence $\Omega^1 A \into P \onto A$ such that
$P$ is a fg projective $A$-bimodule and $\Omega^1 A$ is an
$A$-bimodule which is projective as a left and right $A$-module
(separately), by splicing together
$(\Omega^1 A)^{\o_A j} \into P \o_A (\Omega^1 A)^{\o_A (j-1)} \onto
(\Omega^1 A)^{\o_A (j-1)}$
for all $j \geq 1$.  If we construct these inductively by tensoring on
the left by $\Omega^1 A$, then we see that the sequences are all exact
since $\Omega^1 A$ is a projective right $A$-module; also,
$P \o_A (\Omega^1 A)^{\o_A (j-1)}$ is a projective $A^e$-module
because the result is obvious in the case that $P$ is a free
$A^e$-module.

Now, given $f \in \uHom_{A^e}(\Omega^j A, \Omega^k A)$, we may construct
$\Omega f \in \uHom_{A^e}(\Omega^{j+1} A, \Omega^{k+1} A)$ by applying
$\Omega^1 A \o_A -$, by construction of the above resolution. On the
other hand, the isomorphism $\uHH^i(A, A) \cong \uHH^{i+n}(A, A)$ is
given by the stable module isomorphism $1': \Omega^n A \simeq A$.
That is, we use the stable module isomorphism $\Omega^{n+i} A \simeq \Omega^i A$,
which by the above is 
$1' \o_A \Id: \Omega^n A \o_A \Omega^i A \rightarrow A \o_A \Omega^i
A$, so $\uHH^i(A,A) \cong \uHH^{i+m}(A,A)$ is given by cup product on the left with $1'$.

Now, since cup product with $1'$ induces an isomorphism, we must have
a right inverse $(1')^{-1}$ such that $1' \cup (1')^{-1} = \Id$, and
hence it is a two-sided inverse by Theorem 
\ref{asfhh}.(i),(iii): 
\begin{multline} \label{oneprimes}
  \Id = (1' \cup (1')^{-1}) \cup (1' \cup (1')^{-1}) = 1' \cup
  ((1')^{-1} \cup 1') \cup (1')^{-1} \\ = ((1')^{-1} \cup 1') \cup (1'
  \cup (1')^{-1}) = ((1')^{-1} \cup 1').
\end{multline}
Thus, we also deduce the statement at the beginning of the theorem,
that either $|1'|$ is even, or $2 \cdot \Id = 0$.  (This can also be
deduced from the fact that $2 \cdot (1' \cup 1') = 0$ if $|1'|$ is
odd.)

(ii),(iii) To deduce this, we use that the desired structures exist in
positive degree and satisfy the necessary axioms.  Using the element
$1'$ from part (i) (and Theorem \ref{asfhh}), this result follows
from the following general proposition (see also the comments
after the statement).
\end{proof}
\begin{prop}  \label{calclocp} Let $(\mathcal{V}^\bullet, \mathcal{W}^\bullet)$ be a (pre)calculus and $z \in \mathcal{V}$ a homogeneous element. Then,
there is a unique extension of the calculus structure to the localization
$(\mathcal{V}[z^{-1}]^\bullet, \mathcal{V}[z^{-1}]^\bullet \o_{\mathcal{V}^\bullet} \mathcal{W}^\bullet)$, where by convention, $z^{-1} \wedge z = 1$.
\end{prop}
Now, if $A$ is a periodic Frobenius algebra, with a homogeneous
element $1' \in \uHH^{\geq 1}$ inducing the periodicity, then we claim
that $(\uHH^\bullet(A), \uHH_{\bullet}(A)) = (\uHH^{\geq
  0}(A)[(1')^{-1}], \uHH_{\geq 0}(A)[(1')^{-1}])$.  There is clearly a
map $(\uHH^{\geq 0}(A))[(1')^{-1}] \rightarrow \uHH^\bullet(A)$, which
is an isomorphism in nonnegative degrees, and must therefore be an
isomorphism.  As a result, we deduce that the calculus structure on
$(\uHH^{\geq 0}(A), \uHH_{\geq 0}(A))$ extends uniquely to a calculus
structure on $(\uHH^\bullet(A), \uHH_\bullet(A))$.

We remark that this calculus is not periodic in a trivial way: it is
not true that the Lie derivatives or $B$ must commute with $i_{1'}$.
However, one can write formulas for all the operations in terms of
operations on degrees $0,1,2,\ldots,|1'|-1$ and involving $1'$.

\begin{proof}[Proof of Proposition \ref{calclocp}]
  We use the notation of Definition \ref{calcdfn}, since we are
  discussing general calculi and not only the Hochschild calculus.
    By definition, $\mathcal{V}[z^{-1}]$ is graded commutative.
  Furthermore, note that, since $z \wedge z = (-1)^{|z|} z \wedge z$,
  either $|z|$ is even, or $\mathcal{V}[z^{-1}]$ is an algebra over
  $\Z/2$.  Either way, $z$ is central in $\mathcal{V}[z^{-1}]$ (not
  merely graded-central), and we can omit any mention of $(-1)^{|z|}$.

  Let $\phi: \mathcal{V} \rightarrow \mathcal{V}[z^{-1}]$ denote the
  localization map.  Note that, if $f \in \ker(\phi)$, i.e., $f \wedge
  z^k = 0$ for some $k \geq 0$, then $\{f, g\} \wedge z^{k+1} = 0$ by
  the Leibniz rule, so $\{f, g\} \in \ker(\phi)$. Let $\psi:
  \mathcal{W} \rightarrow \mathcal{W}[z^{-1}] := \mathcal{V}[z^{-1}]
  \o_{\mathcal{V}} \mathcal{W}$ denote the base-change map. Then, we
  have $x \in \ker(\psi)$ iff $\iota_{z^k}(x) = 0$ for some $k \geq
  1$, and similarly to the above, we deduce that $\mathcal{L}_a(x),
  d(x) \in \ker(\psi)$ for all $a \in \mathcal{V}$ using the calculus
  identities. Similarly, for any $y \in \mathcal{W}$, and any $f \in
  \ker(\phi)$, we have $\mathcal{L}_f(y) \in \ker(\psi)$.  Thus, it
  makes sense to speak about the calculus structures as being defined
  on $(\phi(\mathcal{V}), \psi(\mathcal{W}))$, and our goal is to
  extend the structure to $(\mathcal{V}[z^{-1}], \mathcal{W}[z^{-1}])$
  and verify that the calculus identities are still satisfied.

  For operators, we will use $[-,-]$ to denote the \textit{graded}
  commutator: $[\alpha, \beta] := \alpha \circ \beta - (-1)^{|\alpha|
    |\beta|} \beta \circ \alpha$.  For example, $[\LL_a, \iota_b] :=
  \LL_a \circ \iota_b - (-1)^{(|a|-1) |b|} \iota_b \circ \LL_a$.

  We extend the Gerstenhaber bracket from $\phi(\mathcal{V})$ as follows:
\begin{equation}
[z^{-1}, g] := - z^{-2} \wedge [z, g],
\end{equation}
together with skew-symmetry and the Leibniz rule.  We must check that
this yields a well-defined bracket, which amounts to the computation
\begin{equation}
[(f \wedge z) \wedge z^{-1}, g] := z^{-1} \wedge [f \wedge z, g] + (f \wedge z) \wedge [z^{-1}, g] =  [f, g] + z^{-1} \wedge f \wedge [z, g] - (f \wedge z^{-1}) \wedge [z, g] = [f,g].
\end{equation}
It is easy to check that this yields a
Gerstenhaber bracket, and we omit this.

We extend the Lie derivative $\mathcal{L}_-$ to $(\mathcal{V}[z^{-1}]^{\bullet}, \mathcal{W}[z^{-1}]^\bullet)$ as follows.
For $f \in \phi(\mathcal{V})$, we extend the operation $\mathcal{L}_f$ to $\mathcal{W}[z^{-1}]$ using
\eqref{compat1}, with $a := f$ and $b := z$,
using the same procedure as above.  To define the derivative
$\mathcal{L}_{f \wedge z^{-1}}$ we use \eqref{compat2} with $a := f$
and $b := z$.  It is straightforward that this is well-defined.

We must verify that the above gives a precalculus.  We know
that the identities are satisfied when everything is in
$(\phi(\mathcal{V}), \psi(\mathcal{W}))$.  Thus, to verify
that \eqref{compat1} holds, we need to show that the 
LHS (graded)-commutes with
$\iota_{z}$. This follows because $\iota_{b}$ (graded)-commutes with
$\iota_{z}$, and $[\mathcal{L}_a, \iota_{z}] = \iota_{[a,z]}$,
which graded-commutes with $\iota_b$.  We may then inductively 
show that \eqref{compat1} holds:
first, if it holds for any $(a, b)$, it must hold replacing $a$ by
$a \wedge z^{\wedge j}$ for any $j \in \Z$ by construction. Then,
inductively, if the identity holds for $(a, b)$, we may deduce that it
holds for $(a, b \wedge z^{-1})$ using our definition of
$\mathcal{L}_{b \wedge z^{-1}}$ (which is based on \eqref{compat2}).

To verify that \eqref{compat2} holds, it is enough to show that the
identity for $(a, b \wedge z)$ implies the identity for $(a, b)$.  We
may first prove this for $a \in \phi(\mathcal{V})$, and then prove it for all $a$
using the identity for $(z, a)$ and hence for $(a, z)$.  That is,
it suffices prove that the identity for $(a, b \wedge z)$ and $(a, z)$
implies the identity for $(a,b)$. We have
\begin{multline}
\LL_{a \wedge b} \iota_{z} + (-1)^{|a|+|b|} \iota_{a \wedge b} \LL_{z}
= \LL_{a \wedge b \wedge z} 
= \LL_{a \wedge z} \iota_{b} + (-1)^{|a|} \iota_{a \wedge z} \LL_{b} \\ =
\LL_a \iota_{b} \iota_{z} + (-1)^{|a|} \iota_a \LL_{z} \iota_b + (-1)^{|a|} \iota_{a \wedge z} \LL_{b} \\ = \LL_a \iota_b \iota_{z} + (-1)^{|a|+|b|} (\iota_{a \wedge b} \LL_{z} - \iota_a \iota_{[b,z]}) + (-1)^{|a|} (\iota_a \iota_{[z,b]} + \iota_a \LL_b \iota_{z}),
\end{multline}
which gives the identity upon canceling the two inner terms and
multiplying on the right by $\iota_{(z)^{-1}}$.

We can similarly verify that $\LL_{-}$ gives a Lie action. As before,
it suffices
to show that the identity $\LL_{[x,y]} = [\LL_x, \LL_y]$ for
$(a, b \wedge z)$ and $(a, z)$ implies the identity for $(a,b)$.
Since $z$ is invertible, it suffices to verify that $[\LL_a, \LL_b] \iota_{z} = \LL_{[a,b]} \iota_{z}$.   Since (using the Leibniz rule and the fact that
z is central)
\begin{equation}
\LL_{[a, b \wedge z]} + (-1)^{(|a|+1)(|b|+1)}
\LL_{b \wedge [z, a]} = \LL_{[a,b] \wedge z} = 
\LL_{[a,b]} \iota_{z} + (-1)^{|a|+|b|+1} \iota_{[a, b]} \LL_{z},
\end{equation}
it suffices to verify that the LHS equals the RHS after substituting the
desired identity $\LL_{[a,b]} \iota_{z} = [\LL_a, \LL_b] \iota_{z}$.
That is, it suffices to prove that
\begin{equation} \label{deseq}
\LL_{[a, b \wedge z]} + (-1)^{(|a|+1)(|b|+1)}
\LL_{b \wedge [z, a]} = [\LL_a, \LL_b] \iota_{z} + (-1)^{|a|+|b|+1} \iota_{[a, b]} \LL_{z}.
\end{equation}
We have
\begin{gather}
  LHS = [\LL_a, \LL_b \iota_z + (-1)^{|b|} \iota_b \LL_{z}] +
  (-1)^{(|a|+1)(|b|+1)}
  \LL_b \iota_{[z, a]} - (-1)^{|a|(|b|+1)} \iota_b [\LL_{z}, \LL_a], \\
  RHS = [\LL_a, \LL_b] \iota_{z} + (-1)^{|b|} \LL_a \iota_b \LL_{z} - (-1)^{|a| |b|} \iota_b \LL_a \LL_{z},
\end{gather}
from which \eqref{deseq} follows by expanding
$\iota_{[z, a]} = \iota_{z} \LL_a - \LL_a \iota_{z}$
in the first line, and making two pairwise cancellations.

Next, we have to extend differential $d$.
For this, we use the Cartan identity \eqref{cartan}, with $a := z$
(or a power of $z$). We need to check that, with this definition,
\eqref{cartan} holds, and that $d^2 = 0$. We will show that
\eqref{cartan} holds applied to any element $b \in
\mathcal{W}[z^{-1}]$. First, we show that \eqref{cartan} must hold
when $b \in \psi(\mathcal{W})$.  This amounts to the statement that
\begin{equation}
\iota_{[z,b]} = [\iota_{z},[d,\iota_b]].
\end{equation}
We
simplify the RHS as
\begin{equation}
[\iota_{z},[d,\iota_b]] = [[\iota_{z}, d], \iota_b] \mathop{=}^{def.} 
-(-1)^{|z|} [\LL_{z}, \iota_b]
= (-1)^{|z|+|b|(|z|+1)} \iota_{[b,z]} = \iota_{[z,b]},
\end{equation}
as desired (the first step used $[\iota_{z}, \iota_b] = 0$ by
graded-commutativity). Next, we show that \eqref{cartan} must hold for
all $b$.  This amounts to the statement that
\begin{equation}
[d, \iota_{a \wedge z}] = [d, \iota_a] \iota_{z} + (-1)^{|a|} \iota_a
[d, \iota_{z}],
\end{equation}
which follows immediately using $\iota_{a \wedge z} = \iota_a \iota_{z}$.

Finally, to show that $d^2 = 0$, it suffices to show that 
$[d^2, \iota_{z}] = 0$, i.e., $d \LL_{z} + \LL_{z} d = 0$.  Actually,
we know that this identity holds when applied to $\psi(\mathcal{W})$,
so it is enough to show that $[d \LL_{z} + \LL_{z} d, \iota_{z}] = 0$.
We have
\begin{equation}
[d \LL_{z} + \LL_{z} d, \iota_{z}] = -d \iota_{[z,z]} + \LL_{z}^2
+ \LL_{z}^2 - \iota_{[z,z]} d = -\LL_{[z,z]} + 2 \LL_{z}^2 = 0,
\end{equation}
using at the end the fact that $\LL_-$ is an action. 
\end{proof}

Adding the CY Frobenius condition, we obtain the
\begin{thm} \label{pcyfralg}
\begin{enumerate}
\item[(i)]
  Let $A$ be any periodic Calabi-Yau Frobenius algebra. Then, the stable
  Hochschild cohomology is a BV algebra, with BV differential obtained
  by the duality \eqref{hhehh} from the Connes differential. That is,
  the differential $\Delta := \D \circ B  \circ \D$, where $\D$ is the duality \eqref{hhehh} and $B$ the Connes differential, satisfies \eqref{bvid};
\item[(ii)] If $A$ is only $CY(m)$ Frobenius (and not necessarily
  periodic), and $m \geq 1$, then we may still define
  $\Delta := \D \circ B \circ \D$ in degrees $\leq m$, and
  \eqref{bvid} is satisfied when $|a|,|b| \geq 0$ and
  $1 \leq |a|+|b| \leq m$.
\end{enumerate}
\end{thm}
\begin{proof}
(i) The proof is similar to the proof of Theorem 3.4.3 in \cite{Gcy}. Namely,
using \eqref{hhhhint} and \eqref{compat1}, we have
\begin{equation}
\{a, b\} = \D i_{\{a,b\}} \D(1) = a \cup \D \mathcal{L}_b \D(1) -(-1)^{|a|(|b|+1)} \D\mathcal{L}_b\D(a).
\end{equation}
Now, using \eqref{cartan} and \eqref{compat2}, we have
\begin{equation}
RHS = a \cup \Delta(b) - (-1)^{|b|} a \cup b \cup \Delta(1) - (-1)^{|a|(|b|+1)}
 \Delta(b \cup a) + (-1)^{|a| |b|} b \cup \Delta(a).
\end{equation}
Using graded-commutativity, this immediately gives \eqref{bvid}.

(ii) The above proof goes through in the general Calabi-Yau Frobenius case 
in the degrees indicated. Note that we needed $|a|, |b| \geq 0$ because we used \eqref{compat1} applied to $a$ and $b$.
\end{proof}
\begin{rem} \label{bveqcalc} In fact, for any graded-commutative
  algebra $\V^\bullet$, giving the structure of calculus using
  $\W^\bullet := \V^{m-\bullet}$ which satisfies the intertwining property
  \eqref{hhhhint} (for $\D$ the tautological isomorphism) 
  is equivalent to giving $\V^\bullet$ a BV algebra
  structure.  The above theorem showed that calculus+duality gives BV;
  the other direction is as follows: The intertwining property
  \eqref{hhhhint} uniquely specifies what the module structure of
  $\V^\bullet$ on $\W^\bullet$ is, and the differential then gives
  $\LL$.  One may then deduce the remaining identities from
  $\Delta^2 = 0$ and the BV identity \eqref{bvid}.  The identity
  $\LL_{[x,y]} = [\LL_x, \LL_y]$ says precisely that $\Delta^2$ is a
  derivation\footnote{It is worth remarking, by comparison, that, for
    a graded-commutative algebra with an odd differential operator
    $\Delta$ of order $\leq 2$, the Jacobi identity for its principal
    symbol $\pm [\,,]$ says that $\Delta^2$ is a differential
    operator of order $\leq 2$ (so being a derivation is between this
    and $\Delta^2=0$). Skew-symmetry of $[\,,]$ is automatic.};
  \eqref{compat1} says $\Delta$ is a differential operator of order
  $\leq 2$, together with the BV identity (i.e., that
  $a \otimes b \mapsto (-1)^{|a|+1}[a,b]$ is the principal symbol of
  $\Delta$).  Then, \eqref{compat2} is a consequence of
  \eqref{cartan} (and it is a proof that $\Delta$ being an operator of
  order $\leq 2$ yields the Leibniz rule for its principal symbol).
\end{rem}
Put together, any periodic Calabi-Yau Frobenius algebra has Hochschild
cohomology which is a BV algebra and a Frobenius algebra (in the
derived category), and together with Hochschild homology forms a
periodic calculus (together with an isomorphism between the two that
intertwines cup product with contraction). Moreover, the shift functor
$\Omega^{-m}$ is a (right) Serre functor for the category of
fg left modules.

In the case of (generalized, centrally extended) preprojective
algebras (Examples \ref{dynqe},\ref{cedynqe},\ref{gdynqe}) over $\C$
(or any field for the ordinary preprojective algebras), the formula
for the extension of $B$ is quite simple, as we will prove:
\begin{gather} \label{bgsa}
B_\bullet = (-1)^{\bullet}  \D_k \circ B_{-2-\bullet}^* \circ \D_k, \\
B_{-1} = 0. \label{bm10}
\end{gather}
This says that $B$ is \textbf{graded selfadjoint}. As a consequence,
so is the BV differential $\Delta$. This motivates the
\begin{defn} A BV Frobenius algebra is a $\Z$-graded
  Frobenius algebra $H$ ($=$ a $\Z$-graded algebra whose graded
  components are fg projective over $\k$, and with a
  homogeneous invariant inner product of some fixed degree) together
  with a graded self-adjoint differential
  $\Delta_\bullet: H_\bullet \rightarrow H_{\bullet-1}$ of degree
  $-1$, i.e.,
\begin{equation}
\Delta_\bullet = (-1)^{\bullet} \Delta_\bullet^\dagger,
\end{equation}
with $\Delta_\bullet^\dagger$ the adjoint of $\Delta_\bullet$.
\end{defn}
This leads to the
\begin{ques} \label{cyfq1}
Are the formulas \eqref{bgsa},\eqref{bm10} valid for any periodic
Calabi-Yau Frobenius algebra?  
\end{ques}
More generally:
\begin{ques} \label{cyfq2}
  For which Frobenius algebras $A$ does \eqref{bgsa},
  \eqref{bm10} define an extension of the usual calculus to
  $(\uHH^\bullet, \uHH_\bullet)$?
\end{ques}
If the answer to Question \ref{cyfq1} is positive, then the stable
Hochschild cohomology of any periodic CY Frobenius algebra is BV
Frobenius.  If the answer to \eqref{cyfq2} includes CY Frobenius
algebras, then one does not need the periodicity assumption.

To answer the above questions, we suggest to work on the level
of standard Hochschild chains.
Consider the two-sided resolution $\hat N_\bullet$ of $A$ given by
$N_\bullet \onto A \into (N_\bullet^\vee \o_A A^*)$, where $N_\bullet$
is the normalized bar resolution.  Then, the chain complex
$\hat N_\bullet \o_{A^e} A$ has the form
\begin{equation} \label{tsbar}
\cdots \rightarrow A \o \bar A \o \bar A \rightarrow A \o \bar A \rightarrow
A \rightarrow A^* \rightarrow A^* \o (\bar A)^* \rightarrow A^* \o (\bar A)^* \o (\bar A)^* \rightarrow \cdots.
\end{equation}
Then, the two-sided Connes differential $\hat B$ should be given, on chains,
by
\begin{equation} \label{Bhate}
\hat B_i := \begin{cases} B, & i \geq 0, \\
                          0, & i = -1, \\
                          (-1)^{i} B^*, & i \leq -2.
\end{cases}
\end{equation}
We hope to address this in a future paper.

\subsection{Hochschild cohomology of symmetric algebras is BV}\label{symfrobs}
In this subsection, we give a simple
proof that a symmetric Frobenius algebra over an arbitrary base ring
$\k$ has \textbf{ordinary} Hochschild cohomology which is BV.
This is based on using the formulas \eqref{bgsa}, \eqref{bm10} to
extend $B$ (since the algebra need not be periodic, we cannot use 
Theorem \ref{pcyfralg}). 
Actually, we show this more generally for ``stably symmetric''
algebras, i.e., $A \simeq A^\vee$.

First, let us extend the Lie derivative operation to all of $(\uHH^\bullet, \uHH_\bullet)$,
in the spirit of \eqref{bgsa},\eqref{bm10}:
\begin{equation}
\LL_{b} := \hat B i_b - (-1)^{|b|} i_b \hat B.
\end{equation}
We claim that $(i_a \LL_b - (-1)^{|a| (|b|+1)} \LL_b i_a)(x) = i_{\{a,b\}}(x)$,
when $a, b \in \uHH^{\geq 0}(A,A)$ and $x \in \uHH_{\leq -1}(A, A)$. We have
(using graded commutators)
\begin{equation}
  \LL_b(x) = (-1)^{|x| \cdot (|b|+1)} \D_\k (-1)^{|b|} (B^* i_b^*
  - i_b^* B^*)\D_k(x) = -(-1)^{|x| \cdot (|b|+1)} \D_\k [i_b^*, B^*] \D_\k(x).
\end{equation}
The signs above use the identity $c \cdot |x| + d \cdot(|x|+c) =
d \cdot |x| + c \cdot(|x|+d) = (c+d)\cdot |x| + c\cdot d$. (setting $c=1, d = |b|$ or vice-versa). We then obtain
\begin{multline}
[i_a, \LL_b](x) = (-1)^{|x| \cdot (|a|+|b|+1)} \D_\k [[i_b^*, B^*],i_a^*]\D_\k (x)
= (-1)^{|x| \cdot (|a|+|b|+1)} \D_\k ([i_a, [B, i_b]])^*\D_\k (x) \\ = (-1)^{|x| \cdot (|a|+|b|+1)} \D_\k i_{\{a, b\}}^*\D_\k(x) = i_{\{a,b\}}(x).
\end{multline}
But, as in the proof of Theorem \ref{pcyfralg}, this identity
immediately gives the BV identity \eqref{bvid}, letting
$x = \D(1) \in \uHH^{-1}(A, A)$. Note also that, by definition,
$\hat B \D(1) = 0$, so that $\Delta(1) = 0$.

\section{Hochschild (co)homology of ADE preprojective algebras over
  any base}
As mentioned in Example \ref{dynqe}, the preprojective algebra of a
quiver of type ADE is a periodic Calabi-Yau Frobenius algebra.  In
this section, we explicitly describe its Hochschild (co)homology over
an arbitrary base (including positive-characteristic fields). In
characteristic zero, this has already been done in \cite{EE2, Eu}.  We
also review and simplify the algebraic structures on Hochschild
(co)homology, and prove that $\uHH^\bullet$ is BV Frobenius (verifying
Question \ref{cyfq1} in this case). Finally, in \S
\ref{mapsextdyns}, we explicitly describe the maps relating the
Hochschild (co)homology of Dynkin preprojective algebras with the
extended Dynkin case; in the latter case, the Hochschild cohomology
groups were described over $\k = \C$ in \cite{CBEG} and later over
$\Z$ in \cite{S3}.

The new theorems proved here are: Theorem \ref{prepthm}, which
computes the Hochschild (co)homology of the Dynkin preprojective
algebras over $\Z$ and proves that it is a BV Frobenius algebra, and
Theorem \ref{eddcompthm}, which compares the Hochschild cohomology of
the Dynkin and extended Dynkin preprojective algebras over $\Z$.  We
also restate in a simplified form two theorems from \cite{EE2,Eu3}
(Theorems \ref{ee2thm} and \ref{euconnthm}).

We first recall the definition of the preprojective algebras $\Pi_Q$
of ADE type.
\begin{ntn} We will use $[n]$ to denote the degree-$n$-component of a
  graded vector space (discarding other degrees), and $(n)$ to denote
  shifting a graded vector space by degree $n$.  
In particular, $A[m]$ is a graded vector space concentrated in degree $m$,
and $(A(n))[m] =$ $(A[m-n])(n)$. The vector space $A[m](-m)$ is
concentrated in degree zero.
\end{ntn}

\begin{ntn} We will use $\deg(z)$ to denote the degree of an element
  $z$ in a graded algebra or module. This is to distinguish with
  Hochschild degree, where we denote $|a|=m$ if $a \in \uHH^m(A, M)$
  (so $\deg(a)$ would denote the degree with respect to the grading on
  $A$ and $M$).  The notation $[m], (m)$ refer to the $\deg(-)$
  grading, and never to Hochschild degree.
\end{ntn}

Let $\k$ be a commutative ring.  Let $Q$ be a quiver of ADE type with
vertex set $I$.  By convention, $Q$ also denotes the edges of the
quiver.  Let $\dq := Q \sqcup Q^*$ be the double quiver, where
$Q^* := \{a^* \mid a \in Q\}$ is the quiver obtained by reversing all
arrows ($a^*$ is the reverse of $a$).

Let $P_Q, P_{\dq}$ be the associated path algebras over $\k$, and let
$\Pi_Q := P_{\dq} / \ldp r \rdp$ with $r := \sum_{a \in Q} [a,a^*]$.
Let $e_i$ denote the image of the vertex $i$ for any $i \in I$.

Recall from Example \ref{dynqe} that $\Pi_Q$ is a Frobenius algebra over
$\k$ (in fact, periodic Calabi-Yau).
  Let
$(\,,)$ denote an invariant inner product, and let $\eta$ be the
Nakayama automorphism of $\Pi_Q$, so that $(x, y) = (\eta^{-1}(y), x)$.
Recall \cite{RS, ES2, Eu} that we may choose $(\,,), \eta$ such that
\begin{gather}
\eta(e_i) = e_{\bar \eta(i)}, \text{ defining $\bar \eta: I \rightarrow I$ by} \\
\alpha_{\bar \eta(i)} = -w_0(\alpha_i);
\end{gather}
here, $w_0$ is the longest element of the Weyl group of the root system
attached to $Q$, and $\alpha_i, i \in I$ are the roots.  Furthermore,
$\eta$ may be uniquely chosen to act on $\dq \subset P_{\dq}$ so that
\begin{equation}\label{etaform}
\eta(Q^*) \subset \dq, \text{ and } \eta(Q) \subset \begin{cases} - \dq, & \text{if $Q$ is of type $D,E$,} \\  \dq, & \text{if $Q$ is of type $A$},
\end{cases}
\end{equation}
since $Q$ is a
tree.    As a consequence of these formulas, we see that $\eta$ and
$\bar \eta$ are involutions. We remark that $\eta$ is always nontrivial,
even though, for $D_{2n+1}$, $E_7$ and $E_8$, $\bar \eta$ becomes trivial. (Except,
over characteristic $2$, $\eta$ is trivial for $D_{2n+1}, E_7$, and $E_8$.)

Let $m_1, m_2, \ldots, m_{|I|}$ be the exponents of the root system
attached to $Q$, in increasing order. Let $h := m_{|I|}+1$ be the
Coxeter number.

Recall that the Hilbert series of a $\N$-graded vector space $M$ is
defined to be $h(M;t) := \sum_{m \geq 0} \dim M[m] t^m$.  If $E$ is an
$\N$-graded $\k^I$-module for any field $\k$, then we define the
matrix-valued Hilbert series of $E$, $h(E;t)$, by
$h(E;t)_{ij} := \sum_{m \geq 0} \dim E[m]_{i,j} t^d$, where
$E[m]_{i,j} := e_i E[m] e_j$, where $e_i,e_j \in \k^I$ are
the idempotents corresponding to the vertices $i, j \in I$.

\subsection{Reminder of characteristic zero results} \label{ee2res}
Let $A := \Pi_Q$, and assume that $\k$ is a characteristic-zero field. 
We may then describe the Hochschild homology, $HH_*(A)$, and
cohomology, $HH^*(A)$, as follows:
\begin{defn}\cite{EE2} Let
  $U := (HH^0(A)[< h-2])(2), L := HH^0(A)[h-2](-(h-2)), K := HH^2(A)(2)$, and
  $Y := HH^6(A)[-h-2](h+2)$. Also, let $P \in \End(\k^I)$ be the
  permutation matrix corresponding to the involution $\bar \eta$. Let $I^{\bar \eta}$ be the set of vertices fixed by $\bar \eta$.
\end{defn}
\begin{thm}\cite{EE2}\label{ee2thm} (restated)
\begin{enumerate}
\item[(i)] $U$ has Hilbert series
\begin{equation}
h(U;t) = \sum_{i: m_i < \frac{h}{2}} t^{2m_i}.
\end{equation}
\item[(ii)] We have natural isomorphisms $K \cong \ker(P+1)$ and $L \cong \k^{I^{\bar \eta}}$.

\item[(iii)] As graded vector spaces, one has
\begin{gather}
\uHH^0(A,A) \cong U(-2) \oplus Y(h-2), \quad
\uHH^1(A,A) \cong U(-2), \quad
\uHH^2(A,A) \cong K(-2), \\
\uHH^{6+i}(A,A) \cong \uHH^{i}(A,A)(-2h),  \quad \uHH^i(A,A)(2) \cong (\uHH^{5-i}(A,A)(2))^*, \label{prepdual} \\
\uHH_i(A,A) \cong \uHH^{2-i}(A,A)(2), \\
HH^0(A) \cong U(-2) \oplus L(h-2), \quad HH_0(A) \cong \k^I. 
\end{gather}
\end{enumerate}
\end{thm}
\subsubsection{The cup product}
Let us summarize also the cup product structure, which was computed in
\cite{ES2,Eu2}. We explain it using our language and results.  In view of
the first isomorphism of \eqref{prepdual} and Theorem \ref{perfhhthm},
it is enough to consider cup products among elements of Hochschild
degrees between $0$ and $5$.

Since the Calabi-Yau Frobenius dimension is $2$ of shift $2$, by
Theorem \ref{cyfhhthm}, the Hochschild cohomology is a Frobenius
algebra with pairing of Hochschild degree $-5$ (meaning,
$(f, g) \neq 0$ implies that $|f|+|g|=5$ in Hochschild degree), and of
graded degree 4 (meaning, in graded degree, $(f,g) \neq 0$ implies
$\deg(f)+\deg(g) = -4$).  In particular, for all $i \in \Z$, the
composition
\begin{equation}
\uHH^i(A,A) \o \uHH^{5-i}(A,A) \tra^{\cup} \uHH^5(A,A) \tra^{(,\Id)} \k
\end{equation}
is a perfect pairing of graded degree $4$, the same as the second isomorphism
of \eqref{prepdual}.

Moreover, if $|f|+|g|+|h|=5$ (in Hochschild degree), using the
graded-commutativity of cup product, we have
\begin{equation}
  (f \cup g, h) = (-1)^{|g| \cdot |h|} (f \cup h, g) = 
(-1)^{|f| \cdot (|g|+|h|)} (g \cup h, f),
\end{equation}
and since the pairing is perfect, we see that knowing the cup product
$\uHH^{|f|} \o \uHH^{|g|} \rightarrow \uHH^{5-|h|}$ determines the cup
product in the other two pairs of Hochschild degrees, $(|f|,|h|)$ and
$(|g|,|h|)$. That is, we may divide the cup products into the
unordered triples summing to $5$ modulo $6$:
\begin{equation}
(0,0,5), (0,1,4), (0,2,3), (1,1,3), (1,2,2), 
(1,5,5), (2,4,5), (3,3,5), (3,4,4),
\end{equation}
and the cup product in any fixed two Hochschild degrees of a triple
determines the other two pairs of cup products.

The first triple above corresponds to multiplication in the center
$Z(A)$, via the quotient $Z(A) \onto \uHH^0(A, A)$ which performs
$U(-2) \oplus L(h+2) \onto U(-2) \oplus Y(h+2)$ (see \cite{ES2,Eu2}
for an explicit computation of this multiplication). Then, the next
two triples describe $\uHH^1(A, A)$ and $\uHH^2(A, A)$ as
$\uHH^0(A, A)$-modules.  As explained in \cite{ES2,Eu2},
$\uHH^1(A, A)$ is cyclic as an $\uHH^0(A, A)$-module (generated by the
Euler vector field), and since $K$ is concentrated in graded degree
zero, the structure on $\uHH^2(A, A)$ is the obvious (trivial) one: it
is a $\k$-vector space.

The cup products between Hochschild degrees $(1,1)$ and $(3,3)$ are
trivial for graded degree reasons. For types $D,E, A_{2n+1}$, the cup
product is also trivial in Hochschild degrees $(3,4)$ for degree
reasons, and in degrees $(2,4)$ by an argument using the BV identity
(see \cite{Eu2}---the argument only shows that the cup product is
$h$-torsion, and in fact it appears to be nontrivial for type
$D_{2n+1}$ in characteristic two).  In type $A_{2n}$, the cup product
between degrees $(3,4)$ and $(2,4)$ is nontrivial and can be
explicitly computed (see \cite{ES2}; see also the similar type $T$
case in \cite{Eu4}). When nontrivial, the products between degrees
$(3,4)$ and $(2,4)$ are only between the lowest possible degrees (so,
it reduces to a pairing of vector spaces concentrated in bottom
degree, which in fact has rank one since the bottom-degree part in Hochschild degree $4$ has dimension one).

This leaves only the cup products $(1,2,2)$ and $(1,5,5)$. These
are best described as cup products 
\begin{gather}
\uHH^2(A, A) \o \uHH^2(A, A) \rightarrow \uHH^4(A, A) \cong 
\uHH^1(A, A)^*, \\
\uHH^5(A, A) \o \uHH^5(A, A) \rightarrow \uHH^4(A, A) \cong 
\uHH^1(A, A)^*.
\end{gather}
Here, we obtain a nondegenerate symmetric bilinear pairing $\alpha$ on
$K$, and a symplectic pairing $\beta$ on $Y$, respectively (one must
obtain some symmetric and skew-symmetric bilinear pairings on $K$ and
$Y$, respectively, since $K$ and $Y$ are concentrated in degree zero,
and $U$ has $\k$-dimension equal to one in each graded degree; 
nondegeneracy is then a result of explicit computations in
\cite{ES2,Eu2} (Theorem 4.0.8 in \cite{Eu2} for types $D,E$;
throughout \cite{ES2}, part II, for type $A$)).

\subsubsection{The Connes and BV differentials}
Using the dualities and intertwining properties, one immediately
obtains the contraction maps. It remains only to compute the Connes
differential, which yields the BV differential by duality, and then
using the Cartan and BV identities, one immediately computes the Lie
derivatives and Gerstenhaber bracket.

We reprint the Connes differential from \cite{Eu3}. Let
$(z_k), (\omega_k)$ be homogeneous bases for
$U(-2), Y(h-2) \subset \uHH^0(A, A)$, respectively, with $\deg(z_k)= k$.
Let $(\theta_k) \subset U(-2) \subset \uHH^1(A, A)$ be a homogeneous
basis for $\uHH^1(A, A)$ with $\deg(\theta_k) = k$, and let
$(f_k) \subset K \subset \uHH^2(A, A)$ be a basis. We will write
$(f_k^*), (\theta_k^*), (z_k^*), (\omega_k^*)$ for the dual bases of
$\uHH^3(A, A), \uHH^4(A, A)$, and $\uHH^5(A, A)$. (Using dual notation
is where we diverge from \cite{Eu2,Eu3}). Let us abusively identify
these elements with their images under periodicity and $\D$. So, for
example, $\theta_k$ may denote the corresponding element in any group
$\uHH^{1+6s}$ or $\uHH_{1+6s}$ for any $s \in \Z$. (This also differs
from the notation of \cite{Eu2,Eu3}).
\begin{thm}\cite[Theorem 5.0.10]{Eu3}  \label{euconnthm}
The Connes differential $B_\bullet: \HH_\bullet(A,A) \rightarrow \HH_{\bullet+1}(A,A)$ is given by
\begin{gather} \label{euconn1}
0 = B_{6s} = B_{4+6s} = B_{2+6s}(U) = B_{3+6s}(Y^*), \\ \label{euconn2}
B_{1+6s}(\theta_{k}) = \bigl( 1 + \frac{k}{2} + sh)z_{k}, \quad B_{3+6s}(z^*_{k}) = \bigl((s+1)h - 1 - \frac{k}{2}\bigr) \theta^*_{k}, \\ \label{euconn3}
B_{2 + 6s}(\omega_{k}) = \bigl( \frac{1}{2} + s \bigr) 
h \beta^{-1}(\omega_{k}), \\ \label{euconn4}
B_{5+6s}(f_{k}^*) = (s+1) h \alpha^{-1}(f_{k}^*),
\end{gather}
where in the first line, $B_{2+6s}(U)$ means the image of the summand of
$U(2hs)$ under $B_{2+6s}$, and similarly for $B_{3+6s}(Y^*)$.
\end{thm}
Hence, the same formulas are valid for $\uHH_\bullet$ where now $s \in \Z$ is arbitrary.
As a consequence of writing it this way, using the symmetry and
skew-symmetry of $\alpha, \beta$, respectively, it is easy
to verify the (new)
\begin{cor} \label{connsat} The Connes differential $B$ is
  graded selfadjoint with respect to the duality $\Dk$.  Hence, the
  same is true for the BV differential $\Delta$.
\end{cor}
Each line of \eqref{euconn1}--\eqref{euconn4} verifies $B_i = (-1)^i B_{-2-i}^*$ for the concerned summands of $HH_i$. 
\begin{rem}
  In the generalized Dynkin case of type $T_n$ (Example \ref{gdynqe}),
  a similar observation to the above, together with the computation of
  $B$ found in \cite{Eu4}, shows that $B$ is graded selfadjoint in the
  $T_n$ case, and hence so is $\Delta$, i.e., $\uHH^\bullet$ is BV
  Frobenius.
\end{rem}

\subsection{Extension of results to $\Z$ and arbitrary characteristic}
Now, we explain the general $\Z$-structure of Hochschild (co)homology.
Note that, by the Universal Coefficient Theorem, one may immediately
deduce the $\k$-module structure from this for any $\k$; we explain it
for fields $\F_p$ (with $p$ prime) to see the duality. We will also
see that the stable Hochschild cohomology $\uHH^\bullet$ is BV
Frobenius over \textit{any} base field, in Theorem \ref{prepthm}.(v)
below, and give the complete structure of $\Delta$ over any field.
\begin{defn} \label{tukydefn}
  We define (and redefine) the vector spaces $T,U,K,K',Y,Y',L$ by
\begin{gather}
K := \uHH_0(A,A)[0], \quad K' := \tors(\uHH_{-1}(A,A)[0]), \quad T := \uHH_0(A,A)[> 0], \\
Y := (\uHH^{-1}(A,A)[h-2])(-(h-2)), \quad Y' := \tors(\uHH^0(A,A)[h-2])(-(h-2)), \\
U(-2) := \uHH^0(A, A)[<(h-2)], \quad  L := \HH^0(A,A)[h-2](-(h-2)).
\end{gather}
Let $T^*$
  (abusively) denote the graded $\Z$-module
\begin{equation}
T^* := H^1(T^{L*}) \cong \Hom_\Z(T, \Q/\Z),
\end{equation}
where $L*$ denotes the derived dual, and $H^1$ the first cohomology.
\end{defn}
Note that, after tensoring by $\C$, the vector spaces $U, K, Y, L$ above
become (isomorphic to) the ones defined previously.  Also,
$K,K',Y,Y',L$ are all concentrated in degree zero (by construction),
hence the shifts in the definition. On the other hand, $U,T,T^*$ live
in multiple degrees.  Also, note that, by e.g.~the Universal
Coefficient Theorem, when torsion appears in the cohomology over $\Z$,
then if we work instead over $\F_p$, the torsion is tensored by $\F_p$
and replicated in degree one lower (one higher for homology). The free
part carries over (tensored by $\F_p$ of course).
\begin{ntn} Let the \textbf{bad primes} consist of $2$ for $Q=D_n$,
  $2,3$ for $Q=E_6$ or $E_7$, and $2,3,5$ for $Q=E_8$, and no others.
\end{ntn}
\begin{thm} \label{prepthm}
\begin{enumerate}
\item[(i)] 
The modules $T, T^*$ are
finite with torsion elements of order equal to bad primes, and $K',Y'$
are finite with torsion of order dividing the Coxeter number $h$.
The modules $K, Y, L, U$ are free over $\Z$.
\item[(ii)] We have
\begin{gather}
\uHH^0(A, A) \cong U(-2) \oplus Y(h-2) \oplus Y'(h-2), 
\quad \uHH^1(A,A) \cong U(-2), \\
\uHH^2(A, A) \cong K(-2) \oplus T(-2), \quad
\uHH^3(A, A) \cong K(-2) \oplus K'(-2), \\ \uHH^4(A, A) \cong U^*(-2) 
\oplus T^*(-2), \quad \uHH^5(A, A) \cong U^*(-2) \oplus Y(-h-2), \\
\uHH_i(A,A) \cong \uHH^{2-i}(A,A)(2), \quad \uHH^{6+i} \cong \uHH^{i}(-2h), \\
\HH^0(A, A) \cong U(-2) \oplus L(h-2), \quad \HH_0(A, A) \cong \Z^I \oplus
T.
\end{gather}
\item[(iii)] For any prime $p$, letting $M_p := M \o \F_p$ for all $M$, we have
\begin{gather}
\uHH^0(A_p, A_p) \cong U_p(-2) \oplus Y_p(h-2) \oplus Y'_p(h-2), 
\quad \uHH^1(A_p,A_p) \cong U_p(-2) \oplus T_p(-2), \\
\uHH^2(A_p, A_p) \cong K_p(-2) \oplus K'_p(-2) \oplus T_p(-2), \quad \uHH^{5-i}(A_p,A_p)(2) \cong (\uHH_{i}(A_p,A_p)(2))^*, \\
\uHH_i(A_p,A_p) \cong \uHH^{2-i}(A_p,A_p)(2), \quad \uHH^{6+i} \cong \uHH^{i}(-2h), \quad \\ 
\HH^0(A_p, A_p) \cong U_p(-2) \oplus L_p(h-2), \quad \HH_0(A_p, A_p) \cong \F_p^I \oplus
T_p.
\end{gather}
\item[(iv)] For any (bad) prime $p$, $T^*_p$ is cyclic as an $\uHH^0(A_p,A_p)$-module.
\item[(v)] The BV differential $\Delta$ (and hence the Connes
  differential) is graded selfadjoint over any field $\k$. In
  particular, for $\k=\F_p$, $\Delta$ is zero restricted to any
  summand $K'_p, Y'_p, T_p, T^*_p$.
\end{enumerate}
\end{thm} 
The remainder of this subsection is devoted to the proof of the theorem.
\begin{lemma}
The center of $A$ does not increase over positive characteristic. That
is, the map $Z(A) \o \F_p \rightarrow Z(A \o \F_p)$ is an isomorphism
for all ADE quivers.
\end{lemma}
We omit the proof of the above lemma, which was done using MAGMA for
type $E$, and a straightforward explicit computation in the $A_n, D_n$
cases, using bases in terms of paths in the quiver.  Note that the
lemma is actually true for \textbf{all} quivers, since the non-Dynkin
case is proved in \cite{CBEG,S3}.
\begin{cor}
The groups $\HH^1(A, A), \HH^5(A, A), \HH_{1}(A, A),$ and
$\HH_{3}(A, A)$ are all free $\Z$-modules.
\end{cor}
\begin{proof}
  It is enough to show that they are torsion-free. For
  $\HH^1$, this follows from the fact that the differential
  $d_0$ in $C^\bullet(A, A)$ must have saturated image
  (otherwise $\HH^0(A, A)$ would increase in some positive
  characteristic); alternatively, this is a consequence of the
  universal coefficient theorem.  For $\HH^5(A, A)$, we use the
  derived duality
  $\D_\k: C^\bullet(A, A) \iso C^{5-\bullet}(A, A)^*$,
  so that the differential $d_0$ corresponds to $d_4^*$ in the latter
  (again, we could also use the universal coefficient theorem).  Then,
  the duality $\D$ of dimension $2$ and the periodicity $6$ gives the
  results for Hochschild homology: $\HH^j \cong \HH_{6n+2-j}$ whenever
  $j, 6n+2-j$ are both positive.
\end{proof}
Using the duality $\D$ and the periodicity by period $6$, it
remains only to compute the torsion of
$\HH^2(A, A), \HH^3(A, A), \HH^4(A, A)$, and
$\HH^6(A, A)$.  Using the duality $\D_\k$, $\HH^2$ and $\HH^4$
must have dual torsion, so it is really enough to compute in degrees
$2, 3,$ and $6$.  We will see that the torsion in these degrees
will be nontrivial, but only  in bad
primes ($2$ for $D_n$, $2,3$ for $E_6$ and $E_7$, and $2,3,5$ for
$E_8$), and primes dividing the Coxeter number $h$.

First, by the duality $\D$, the torsion of $\uHH_0(A, A)$ and
$\uHH^2(A, A)$ are isomorphic. Since $\HH_0(A, A)$ has no torsion in
degree zero, and the inclusion $\uHH_0(A, A) \subset \HH_0(A, A)$ is
full in nonzero degrees, the torsion of $\uHH_0(A, A)$ (and hence of
$\uHH_2(A, A)$) is the same as that of $\HH_0(A, A)$.  The latter was
computed in \cite{S3}, and we collect results for convenience:
\begin{prop}\cite[Theorem 4.2.60]{S3} \label{s3pr} 
The module $\HH_0(A, A) \cong \Z^I \oplus T$, where
$T := \HH_0(A, A)_+$ is finite and given
as follows:
\begin{itemize}
\item For $Q = A_n$, $T = 0$, 
\item For $Q = D_n$, 
\begin{equation}
 T \cong \bigoplus_{4 \mid m, 0 < m \leq 2(n-2)} \Z/2(m),
\end{equation} 
$[r_2 \cdot z^{m-1}]$ for $\deg(r_2)=\deg(z)=4$, where
\item For $Q = E_n$, $T$ is a (finite) direct sum 
of shifted copies of $\Z/2$ and $\Z/3$, and in the case $n=8$, 
also of $\Z/5$. In particular:
\begin{gather}
T_{E_6}  \cong \Z/2(4) \oplus \Z/3(6), \\
T_{E_7} \cong (\Z/2(4) \oplus \Z/2(8) \oplus \Z/2(16)) \oplus \Z/3(6), \\
T_{E_8} \cong (\Z/2(4) \oplus \Z/2(8) \oplus \Z/2(16) 
\oplus \Z/2(28)) \oplus (\Z/3(6) \oplus \Z/3(18)) \oplus \Z/5(10).
\end{gather}
\end{itemize}
Moreover, for any $Q$ and any bad prime $p$, there exists a top-degree
torsion element, $\rptop$, such that all homogeneous $p$-torsion
elements $[x]$ have the property that $[x \cdot z] = \rptop$ for
some homogeneous central element $z \in \HH^0(\Pi_Q, \Pi_Q)$.
\end{prop}
We deduce immediately that $T^* \o \F_p$ is cyclic under dual contraction,
and hence (by \eqref{dkint}) also under contraction. By the intertwining
property \eqref{hhhhint}, we deduce part (iv) of Theorem \ref{prepthm}.

It remains only to compute the torsion of $\uHH^3(A, A)$ and
$\uHH^6(A, A)$, i.e., to compute $K'$ and $Y'$ and verify that there
is no other torsion. We will use some results of \cite{Eu2} for this, but
let us explain them using our language.  
Using the formulation of Definition \ref{tukydefn}
and the (normalized) bar complex, it suffices to compute the
torsion of the cokernels of the maps
\begin{equation} \label{cokmapseq}
C_{0}(A,A) \cong A \tra^{d_0}  A^* \cong C_{-1}(A,A), \quad C^{-1}(A,A) 
\cong A_{\eta^{-1}} \tra^{d^{-1}} A \cong C^0(A,A).
\end{equation}
To express these maps, write $\Id = \sum_i x_i^* \o x_i \in A^* \o A$.
Then, the map $A \o A \rightarrow A^* \o A$ in the normalized bar
resolution is given by $(x \o y) \mapsto xy \cdot \Id$.  So, the
``conorm'' differential $d_0$ in \eqref{cokmapseq} must be given by
$y \mapsto \sum_i x_i y x_i^* \in A^*$.  Since we may assume
$\deg(x_i) = -\deg(x_i^*)$, the image can only be nonzero if
$\deg(y)=0$.  Similarly, the ``norm'' differential $d^{-1}$ is given
by $d^{-1}(y) = \sum_i x_i^* y x_i \in A$, now viewing $x_i^*$ as an
element of $A$ via $(\,, )$.  Here also, only $\deg(y)=0$ need be
considered. We deduce
\begin{prop}
\begin{enumerate}
\item[(i)]
The conorm, $d_0$, and norm, $d^{-1}$, maps are given by
\begin{gather} \label{nocon}
d_0(y) = \sum_i x_i y x_i^* \in A^*, \quad d^{-1}(y) = \sum_i x_i^* y x_i \in A,\intertext{where}  \Id = \sum_i x_i^* \o x_i \in A^* \o A.
\end{gather}
The image of the conorm and norm maps $d_0, d^{-1}$ must lie in the
top degree of $A^*, A$, respectively.  
\item[(ii)] Identifying $A^*$ with $A$ using $(\,,)$, and letting
$\omega_i \in e_i A e_{\bar \eta(i)}$ denote the image of $e_i^* \in A^*$,
we have
\begin{equation} \label{normcoflas}
d_0(e_i) = \sum_{j \in I} \tr(\Id|_{e_j A e_i}) \omega_j, \quad
d^{-1}(e_i) = \sum_{j \in I}\delta_{j,\bar \eta(j)} \delta_{i, \bar \eta(i)} 
\tr(\eta|_{e_i A e_j}) \omega_j.
\end{equation}
\end{enumerate}
\end{prop}
\begin{proof} Part (i) has already been proved and is more generally
  true for any (graded) Frobenius algebra $A$. We show part (ii), 
  using the fact that $A[0] \cong \k^I$ as a subalgebra, with $\eta$ acting
  by the permutation $\bar \eta$.  Let
  $f: A \rightarrow \k$ be the function such that $(a,b) = f(ab)$.  We
  have that $f(e_i, \omega_i) = f(\omega_i) = 1$ for all $i$.  Thus,
  to compute $d_0(e_i)$, it is enough to find
  $f(e_j d_0(e_i))= f(\sum_{\ell} (e_j x_\ell e_i x_\ell^*))$ for all
  $j \in I$.  This is the same as
  $\sum_{\ell} (e_j x_\ell e_i, x_\ell^*)$, which is the trace of
  the projection $A \onto e_j A e_i$.  Similarly, we have
\begin{equation}
\sum_\ell (e_j x_\ell^* e_i, x_\ell) = \sum_\ell (\eta(e_i x_\ell) e_j, x_\ell^*) = \tr(x \mapsto e_{\bar \eta(i)} \eta(x) e_j) = \tr(\eta|_{e_i A e_j}) 
\delta_{j,\bar \eta(j)} \delta_{i, \bar \eta(i)}.
\end{equation}
\end{proof}
We note that in \cite{Eu2}, sums such as \eqref{nocon} (with $x_i$ a
basis) are used to describe $K$ and $Y$; the proposition above
explains their origin through norm and conorm maps. In particular,
bases are not needed, and under the connectivity assumption
$A[0]\cong \k^I$ of (ii), one can re-express the sum as a trace. We
believe that the necessity of using such formulas to describe usual
Hochschild cohomology gives further justification for studying stable
Hochschild cohomology.

Using the proposition, to compute $K'$ and $Y'$,
it suffices to compute two matrices: an $I \times I$-matrix,
$H_Q$, whose entries are $(H_Q)_{ij} = \tr(\Id_{e_i A e_j})$, and
an $I^{\bar \eta} \times I^{\bar \eta}$-matrix, $H_Q^\eta$, whose
entries are $(H_Q^\eta)_{ij} = \tr(\eta_{e_i A e_j})$, where
$I^{\bar \eta} := \{i \in I: \bar \eta(i)=i\}$.  These matrices were computed
in the $D,E$ cases in \cite{Eu2}.  In fact, $H_Q$ itself was originally
computed for all Dynkin cases in \cite{MOV}: since $A$ is free over $\k$,
$H_Q = h(A;1)$ is the Hilbert series matrix $h(A;t)$ evaluated at $t=1$. In \cite{MOV} is the following formula for $h(A;t)$:
\begin{equation}
h(A;t) = (1 + P t^h)(1 - Ct + t^2)^{-1}.
\end{equation}
So $H_Q = (1+P)(2-C)^{-1}$.  These are then easy to compute.
It is also not difficult to compute $H_Q^\eta$, which we omit, since
the $D,E$ cases are already in \cite{Eu2}, and the $A$ case is easy.
We obtain the following (for $\k=\Z$):
\begin{prop} 
$K'$ and $Y'$ are zero if $Q = A_n$, and otherwise are given by 
\begin{gather}
K' \cong \begin{cases}  (\Z/2)^{\oplus 2\lfloor \frac{n}{2} \rfloor -2}, & Q = D_n, \\
                       (\Z/2)^{\oplus 2}, & Q = E_6, \\
                       (\Z/2)^{\oplus 6}, & Q = E_7, \\
                       (\Z/2)^{\oplus 8}, &  Q = E_8.
\end{cases} \\
Y' \cong \begin{cases} \Z/2, & Q=D_n, \quad n \text{ even}, \\
                       \Z/(n-1)=\Z/(h/2), & Q = D_n, \quad n \text{ odd}, \\
                       \Z/3, & Q = E_7, \\
                       0, & \text{ otherwise.}
\end{cases}
\end{gather}
\end{prop}
This proves (i) of Theorem \ref{prepthm}. At this point, (ii) and
(iii) are immediate from the dualities and the Universal Coefficient Theorem.

It remains to prove part (v).  Since we already know (by Corollary
\ref{connsat}) the result for characteristic zero, it suffices to take
$\k=\F_p$ for some prime $p$.  As before, let $1 \in \uHH^0$ be the
identity, and $1' \in \uHH^6$ induce the periodicity. First, we note
that $\Delta(1) = 0$ and $\Delta(1') = 0$ because this is true over
$\Z$ (by \cite{Eu3}), and moreover that $\{(1')^i, (1')^j\} = 0$ for
all $i,j$ for the same reason.  Hence, $L_{1'}$ is also graded
self-adjoint, and it suffices to verify that
$(\Delta(a), b) = (-1)^{|a|} (a, \Delta(b))$, when
$0 \leq |a|, |b| \leq 5$, and hence either $|a|+|b| = 6$ or
$|a|=|b|=0$.  To do this, we will show that $\Delta$ kills summands of
the form $K'_p, Y'_p, T_p$ and $T_p^*$ (the second statement of (v)).

It is clear, for graded degree reasons, that $\Delta$ kills summands of the
form $Y', K', T, T^*$ over $k = \Z$. It remains to show that the new
summands appearing over $\F_p$ are also killed. For $K'$ in $\uHH^2$,
this is true for degree reasons and the fact that the kernel of
$\Delta$ on $\uHH^1[0]$ is zero (from the characteristic-zero case),
using that $\Delta^2 = 0$. For degree reasons, the summand of $T_p^*(-2)$
in $\uHH^3$ and the summand of $Y'(-h-2)$ in $\uHH^5$ must be killed.

It remains only to show that the summand of $T_p(-2)$ in $\uHH^1$ is killed.
Note first that there is an element
$HEu \in \uHH^1(\Pi_Q)$, the ``half-Euler vector field'', whose action
on closed paths (which must have even length since $Q$ is a tree) is
to multiply by half the path-length.  It follows from \cite{Eu2}
(cf.~\cite{S}, \S 10 for the extended Dynkin case and the following
subsection) that $HEu$ generates $\uHH^1(\Pi_Q)$ as a
$\uHH^0(\Pi_Q)$-module in the case that $\k = \Z$. Now, over $\k=\F_p$,
$U(-2) \o \F_p$ is isomorphic to a direct summand $U_p(-2)$ of $\uHH^1(\Pi_Q)$.
From the argument in \cite{Eu3}, we know that the operator $\Delta|_{\uHH^1(A,A)}$ acts on $U_p(-2)$ by $\Delta(z HEu) = (\frac{\deg(z)}{2}+1) z$, for all $z \in \uHH^0(\Pi_Q)[< h-2]$.

We claim that the Lie derivative $L_{HEu}$ acts on
$\uHH_\bullet$ by multiplication by half the graded degree.  For
$\HH_\bullet$, this follows from the explicit formula for the Lie
derivative (as argued in \cite{Eu3}), and then this extends easily to
$\uHH_\bullet$ by taking the unique extension guaranteed by
Proposition \ref{calclocp} (using the construction given in the
proof).  It follows from this that, under the Gerstenhaber bracket,
$\text{ad}(HEu)$ acts on $\uHH^\bullet$ by multiplication by half the graded
degree.

From the identity \eqref{bvid} and the fact that $\Delta(1) = 0$, we have
\begin{equation}
\Delta(HEu \cup x) - \Delta(HEu) \cup x + HEu \cup \Delta(x) =  \{HEu, x\} = \frac{\deg(x)}{2} x,
\end{equation}
for any homogeneous $x \in \uHH^1$. Since $\Delta(HEu)=1$ and $\Delta|_{\uHH^2} = 0$, we obtain
\begin{equation}
HEu \cup \Delta(x) = (\frac{\deg(x)}{2}+1) x.
\end{equation}
For any choice of splitting of $T_p$ in $\uHH^1$, the structure of $T$
as given in Proposition \ref{s3pr} shows that
$\frac{\deg(x)}{2}+1 = 0$ as an element of $\F_p$, for all $x \in
T_p$. Hence, we have
\begin{equation}
HEu \cup \Delta(x) = 0.
\end{equation}
Since, for degree reasons, $\Delta(x)$ is in the $U_p(-2)$-summand of
$\uHH^0$, we deduce that $\Delta(x) = 0$.  
This completes the proof of Theorem \ref{prepthm}.

\subsection{The maps between the extended Dynkin and Dynkin 
preprojective algebras} \label{mapsextdyns}
In this subsection, we will interpret $K', Y',$ and $T$ in terms of
the preprojective algebras of the corresponding extended Dynkin
quivers.  For this, we use the projection $\pi: \Pi_{\tilde Q} \onto \Pi_Q$
and functoriality of Hochschild homology.  We need to recall a few
results from \cite{S3} first.

Let $Q$ be an ADE quiver and let $\tilde Q$ be the corresponding
extended Dynkin quiver.  By \cite{CBEG,S3}, we know that
$\Pi_{\tilde Q}$ is (ordinary) Calabi-Yau of dimension $2$, and in
particular has Hochschild dimension $2$.  Let
$Z_{\tilde Q} := \HH^0(\Pi_{\tilde Q}, \Pi_{\tilde Q})$ and
$Z_{\tilde Q, +} := Z_{\tilde Q}[\geq 1]$. The ring $Z_{\tilde Q}$ is
closely related to the Kleinian singularity ring: one has
$Z_{\tilde Q} \o \C \cong \C[x,y]^\Gamma$ where $\Gamma$ is the group
corresponding to $\tilde Q$ under the McKay correspondence, and one
can even replace $\C$ with
$\Z[\frac{1}{|\Gamma|}, e^{\frac{2 \pi i}{|\Gamma|}}]$.  There is a
standard integral presentation of $\C[x,y]^\Gamma$ which actually
describes $Z_{\tilde Q}$ over $\Z$ (see, e.g., \cite[Propositions 6.4.2, 7.4.1, and 8.4.1]{S3}).

By \cite[\S 10.1]{S3}, we know that, as graded $\Z$-modules,
\begin{gather}
\HH_0(\Pi_{\tilde Q},\Pi_{\tilde Q}) \cong \Z^I \oplus Z_{\tilde Q, +} \oplus T, \\
\HH_1(\Pi_{\tilde Q}, \Pi_{\tilde Q}) \cong Z_{\tilde Q}(2) \oplus  Z_{\tilde Q, +} \\
\HH_2(\Pi_{\tilde Q}, \Pi_{\tilde Q}) \cong Z_{\tilde Q}(2), \quad \HH^i \cong \HH_{2-i}.
\end{gather}
\begin{thm} \label{eddcompthm}
The induced maps $\pi_{*, i}: \HH_i(\Pi_{\tilde Q}, \Pi_{\tilde Q}) \rightarrow 
\HH_i(\Pi_Q, \Pi_Q)$ are given as follows:
\begin{itemize}
\item[(0)] $\pi_{*,0}$ is an isomorphism on $\Z^I \oplus T$, and kills 
$Z_{\tilde Q,+}$;
\item[(1)] $\pi_{*,1}$ is a surjection $Z_{\tilde Q}(2) \onto U$, with kernel
the elements of degree $\geq h$ (and killing the second factor, $Z_{\tilde Q,+}$);
\item[(2)] $\pi_{*,2}: Z_{\tilde Q}(2) \rightarrow U \oplus Y(h) \oplus Y'(h)$
is a surjection onto $U \oplus Y'(h)$, killing $Z_{\tilde Q}[> h-2]$, and
sending $Z_{\tilde Q}[h-2]$ onto torsion.
\end{itemize}
Moreover, these maps give rise to maps
$\HH^i(\Pi_{\tilde Q}, \Pi_{\tilde Q}) \rightarrow \HH^i(\Pi_Q,
\Pi_Q)$
for $i \in \{0,1\}$, which describe the image of central elements, and
describe descent of outer derivations, related to the above by $\D$.
On $\HH^0$, the map $Z_{\tilde Q}[h-2] \rightarrow L(h-2)$ maps
isomorphically to the saturation of the kernel of
$L(h-2) \onto (Y(h-2) \oplus Y'(h-2))$ (i.e., the kernel of
$L(h-2) \onto Y(h-2)$).
\end{thm}
\begin{proof}
Part (0) is a consequence of \cite[Theorem 4.2.60]{S3}.

(1) Let us prove this by computing instead the map on $\HH^1$; there
are duality isomorphisms $\HH^1 \iso \HH_1$ which must commute with
$\pi_*$ because they are realized by explicit maps of chain complexes
expressed in terms of the quiver.  Recall also that $\HH^1(A, A)$ is
the space of outer derivations (derivations of $A$ modulo inner
derivations).

We claim that the outer derivations descend from $\Pi_{\tilde Q}$ to
$\Pi_Q$. To see this, we use the explicit description of them for
$\Pi_{\tilde Q}$ from \cite[\S 10.2]{S3}: the outer derivations are
realized by certain $\Q$-linear combinations of the half-Euler vector
field (multiplying by the degree in $Q$, setting degree in $Q^*$ to be
zero), and maps $\phi_x: y \mapsto \{x,y\}$, using the Poisson bracket
$\{\,,\}$ induced by the necklace Lie bracket on $\HH_0$ (or the
Poisson structure on $\C[x,y]^\Gamma$). The latter was shown to make
sense as a map
$\HH_0(\Pi_{Q'}, \Pi_{Q'}) \o \Pi_{Q'} \rightarrow \Pi_{Q'}$ in
\cite[\S 5.2]{S3}, for any quiver $Q'$, and the half-Euler vector
field obviously makes sense. Also, although $\HH^1$ consists of some
fractions of sums of these outer derivations, clearly an outer
derivation is a multiple of some integer on $\Pi_{\tilde Q}$ only if
the same is true in $\Pi_Q$.
 
Next, we claim that all outer derivations on $\Pi_Q$ are obtained in this
  way. This is an immediate consequence of the fact  that
  $\HH^0(\Pi_{\tilde Q}, \Pi_{\tilde Q})$ surjects to $U(-2)$, since
  $\HH^1(\Pi_Q, \Pi_Q)$ is $U(-2)$ times the half-Euler vector
  field mentioned above (cf.~\cite{Eu2}, Proposition 8.0.4, and
  \cite[II]{ES2}).  We thus deduce the desired statement, and (1).

  Next, we prove (2). We note that this is equivalent to the desired
  statement on the level of $\HH^0$, i.e., for the map
  $Z_{\tilde Q} \rightarrow Z_{Q}$, by virtue of the duality maps
  $\D: \HH^0(\Pi_{\tilde Q}, \Pi_{\tilde Q}) \iso \HH_2(\Pi_{\tilde
    Q},
  \Pi_{\tilde Q})$
  and $\D: \HH^0(\Pi_Q, \Pi_Q) \onto \HH_2(\Pi_Q, \Pi_Q)$, where the
  latter is the quotient $L \onto (Y \oplus Y')$ on $L$, and the
identity on $U$.

  To prove (2), we use the fact that the Connes differential, and
  hence the BV differential, are functorial.  For the extended
  Dynkin side, by \cite[Theorem 10.3.1]{S3}, the BV differential
  $\Delta: \HH^1 \rightarrow \HH^0$ is the map sending
  $z \cdot HEu$ to $(HEu+1)z$, where $HEu$ is the half-Euler vector
  field, and $z \in HH^1$; $\Delta$ kills the derivations related to
  the Poisson bracket as above. In other words, the map
  $B: \HH_1 \rightarrow HH_2$ sends
  $(z, w) \in Z_{\tilde Q}(2) \oplus (Z_{\tilde Q})_+$ to
  $(HEu+1)(z) \in Z_{\tilde Q}$ (here we ignored the shift by two in
  applying $HEu$).  On the Dynkin side, the Connes differential is also
  given by $z \mapsto (HEu+1)(z)$, for $z \in U$.  Hence, by
  functoriality of the Connes differential, we deduce that $\pi_{*,2}$
  is as described in degrees $< h$, and in degree $h$, has to at least
  map to $Y'$ (there can be $h/2$-torsion on account of the $(HEu+1)$).

  To complete the argument, it suffices to prove the surjectivity to $Y'$:
  in terms of $\HH^0$, we have to show that the map
  $Z_{\tilde Q}[h-2] \rightarrow L(h-2)$ maps surjectively to the
  kernel of $L(h-2) \onto Y(h-2)$.  For this, we can perform a
  relatively easy explicit computation, showing that the elements from
  \cite{S3} map to the saturation of the column span of $H_Q^\eta$.
  For $A_n, D_n$ this is straightforward; for $E_6, E_8$, there is
  nothing to show; and for $E_7$, where $h-2=16$, this alternatively
  follows from Proposition 7.3.3 of \cite{Eu} (which computes the square of
  an element $z_8 \in Z_{\tilde Q}[8]$: this turns out to be the
  needed element which spans the rank-one kernel of $L \onto Y$. But,
   we already know
  that $Z_{\tilde Q}[8] \iso Z_Q[8]$ by the above.)
\end{proof}
\begin{rem}
  The above gives an alternative (integral) computation of the algebra
  structure on $\HH^0(\Pi_Q, \Pi_Q)$ given in \cite{Eu2}, \S 7: this
  must be obtained from truncating the ``Kleinian singularity''
  algebra $\HH^0(\Pi_{\tilde Q}, \Pi_{\tilde Q})$ at degrees
  $\leq h-2$, and composing with the kernel map
  $\ker(L \onto Y) \into L$. The asserted relation to the Kleinian
  singularity $\C^2/\Gamma$ associated to $Q$ is that
  $Z_{\tilde Q} \cong e_{i_0} \Pi_{\tilde Q} e_{i_0}$ (\cite{S},
  Theorem 10.1.1), where $i_0$ is the extending vertex of $\tilde Q$,
  and that $e_{i_0} \Pi_{\tilde Q} e_{i_0} \o \C \cong \C[x,y]^\Gamma$
  (alternatively, instead of by $\C$, one can tensor by any ring
  containing $\frac{1}{|\Gamma|}$ and $|\Gamma|$-th roots of unity).
\end{rem}
\begin{rem}
  In \cite{S3}, the exact structure of $T$ turned out to be mandated
  by the requirement that, for non-Dynkin, non-extended Dynkin quivers
  $\hat Q \supsetneq \tilde Q \supsetneq Q$, the torsion of
  $\HH_0(\Pi_{\hat Q}, \Pi_{\hat Q})$ is $\Z/p$ in degrees $2p^\ell$
  for all primes $p$ and all $\ell \geq 1$, and these are generated by
  elements of the form $\frac{1}{p}[r^{p^\ell}]$ (where
  $r = \sum_{a \in Q} a a^* - a^* a$ is the relation).  The specific
  structure of the torsion in the Dynkin and extended Dynkin cases
  compensates for the fact that $Z_{\tilde Q}$ is missing some degrees
  that would otherwise be necessary to produce the torsion of
  $\HH_0(\Pi_{\hat Q}, \Pi_{\hat Q})$ (using the description
  in \cite[Theorem 4.2.30]{S3},  of torsion
  elements of $\HH_0(\Pi_{\hat Q}, \Pi_{\hat Q})$ not coming from
  $\HH_0(\Pi_Q, \Pi_Q)$ as cyclic products of elements
  $r_{\tilde Q}$ with elements of $Z_{\tilde Q}$.)
\end{rem}

\section{Hochschild (co)homology of centrally extended preprojective algebras} 
\label{ces}
In this section, we compute the BV structure on the Hochschild
cohomology of centrally extended preprojective algebras $A$ over $\k =
\C$, and verify that the BV differential is graded self-adjoint
(hence, the Hochschild cohomology is a BV Frobenius algebra).  From
this, the structure of calculus on $(\uHH^\bullet(A),
\uHH_\bullet(A))$ easily follows as in Remark \ref{bveqcalc}, using
the duality $\D$ (we omit the explicit formulas).

As before, let $Q$ be a quiver of $ADE$ type with vertex set $I$.  In
\cite{ER}, the centrally extended preprojective algebra $\Pi_Q^\mu$ is
defined as a central extension of $\Pi_Q$, in terms of a parameter
$\mu \in \k^I$.  We assume that $\mu$ is a \textbf{regular weight},
i.e., if $\mu = \sum_{i \in I} \mu_i \cdot e_i$ for
$\{e_i\} \subset \k^I$ the idempotents corresponding to $I$, then
$(\sum_i \mu_i \omega_i, \alpha) \neq 0$ for any root $\alpha$ of the
root system attached to $Q$, with $\omega_i$ the fundamental weights.
Explicitly, for all $\alpha = \sum_{i \in I} \alpha_i e_i \in \Z^I$
such that
$\sum_{i \in I} \alpha_i^2 - \sum_{a \in Q} \alpha_{h(a)}
\alpha_{t(a)} = 1$, we have $\sum_i \alpha_i \mu_i \neq 0$.

Let $P_\dq[z]$ be the algebra of polynomials in the
central parameter $z$ with coefficients in $P_\dq$. Then we define
\begin{equation}
\Pi_Q^\mu := P_\dq[z] / \langle z \cdot \mu - \sum_{a \in Q} [a, a^*]\rangle.
\end{equation}
This is a graded algebra with $\deg(e) = \deg(e^*) = 1, \deg(z) = 2$, for all
$e \in Q$. Now, let $A := \Pi_Q^\mu$, and let $Z$ be the center of
$A$.  Let $h$ denote the Coxeter number of $Q$. Let $A_+ := A[\geq 1]$
denote the part of positive degree, and let $A_{top}:= A[2h-2]$ denote
the part of $A$ of top degree.

In \cite{ER}, it is proved that $\Pi^\mu$ is Frobenius over $\k=\C$.
We note that, over $\k=\Z$, this is not, in general, true.  For
example, for $Q=A_2$ and $\mu = \rho = \sum_{i \in I} e_i$, we obtain
that $\Pi_{A_2}^\rho \cong P_\dq / 2 P_\dq[\geq 3]$, which is not fg
or projective over $\Z$.  However, certain parameters $\mu \in \Z^I$
should yield a Frobenius algebra, and we hope to explore this in a
future paper.  Namely, these parameters should be those such that
$(\sum_{i \in I} \mu_i \omega_i, \alpha) = \pm 1$ for all roots
$\alpha$; more generally, for $\Pi_{Q}^\mu \otimes \F_p$ to be
Frobenius over $p$, the condition should be that
$(\sum_{i \in I} \mu_i \omega_i, \alpha)$ is not a multiple of $p$.
We hope to explore this in a future paper.

For the rest of this section, let us take $\k := \C$ and assume that
$\mu$ is regular ($(\sum_i \mu_i \omega_i, \alpha) \neq 0$).

There is a periodic resolution of $A$ of period $4$ (\cite{Eu}, \S 3),
and $A$ is a symmetric algebra, so we immediately deduce (as stated in
Example \ref{cedynqe}) that $A$ is a periodic Calabi-Yau Frobenius
algebra of dimension $3$ (of shift $4$) and period $4$ (of shift
$2h$).
\begin{thm}\cite{Eu} The Hochschild cohomology groups of 
$A$ over $\k := \C$ are given by (for $n \geq 0$):
\begin{gather}
HH^{4n+1}(A) \cong (Z \cap \mu^{-1}[A,A])(-2nh-2) \cong z Z(-2nh-2), \\
HH^{4n+2}(A) \cong A/([A,A] + \mu Z)(-2nh-2), \\
HH^{4n+3}(A) \cong A_+/[A,A](-2nh-4), \\
HH^{4n+4}(A) \cong Z/A_{top}(-2(n+1) h). 
\end{gather}
\end{thm}
From the periodicity, we immediately deduce the groups $\uHH^\bullet$,
by allowing $n$ to be an arbitrary integer in the above. The fact that
$\uHH^i \cong (\uHH^{3-i})^*$ says that the nondegenerate trace pairing
\cite{ELR} induces nondegenerate pairings
\begin{gather}
(Z \cap \mu^{-1}[A,A])(-2) \o A/([A,A] + \mu Z) \rightarrow \k, \\
A_+/[A,A](-4) \o Z/A_{top} \rightarrow \k.
\end{gather}
To describe the cup products, as before, it suffices to describe the
product between two degrees for every triple of integers between $0$
and $3$, which sums to $3$ modulo $4$:
\begin{equation} \label{cecuppos}
(0,0,3), (0,1,2), (1,1,1), (1,3,3), (2,2,3).
\end{equation}
Recall from \cite{ELR,Eu} the Hilbert series for these graded vector
spaces, using
again $m_1 < \ldots < m_{|I|} = h$ to denote the exponents of the root
system (note that the sets $\{h-m_i\} = \{m_i\}$ are identical):
\begin{gather} \label{hsce1}
h(\uHH^0(A);t) = h(\uHH^1(A);t) = 
\sum_{i=1}^r \bigl( t^{2m_i-2} + t^{2m_i} + \cdots + t^{2h-6} \bigr), \\ \label{hsce2}
h(\uHH^2(A);t) = h(\uHH^3(A);t) = 
\sum_{i=1}^r \bigl( t^{-2} + 1 + \cdots + t^{2m_i - 6} \bigr).
\end{gather}
\begin{thm} \label{cprepthm}
\begin{enumerate}
\item[(i)] As modules over $\uHH^0(A)$, we have
  $\uHH^1(A) \cong \uHH^0(A)$ and
  $\uHH^2(A) \cong \uHH^3(A) \cong (\uHH^0(A))^*$.
\item[(ii)] All of the cup products
  $\uHH^i(A) \o \uHH^j(A) \rightarrow \uHH^{i+j}(A)$ for
  $1 \leq i \leq j \leq 3$ are zero except for
  $\uHH^1(A) \o \uHH^2(A) \rightarrow \uHH^3(A)$, which, using the
  identifications of (i), is the canonical map
  $\uHH^0(A) \o (\uHH^0(A))^* \rightarrow (\uHH^0(A))^*$.
\end{enumerate}
\end{thm}
\begin{proof} (i) From \cite{Eu}, p. 10, it follows easily that
  $\uHH^1(A)$ is of the desired form.  Since there must be a unique
  (up to scaling) element in degree zero, we can use the explicit
  isomorphism $\uHH^0(A) \iso \uHH^1(A)$, $z \mapsto z \cdot Eu$ where
  $Eu$ is the Euler vector field.  Then, the statements about
  $\uHH^2(A), \uHH^3(A)$ follow immediately from Theorem
  \ref{cyfhhthm}.  (We may even show compatibility with the duality
  pairings defined in \cite{Eu} using the trace map of \cite{ELR}, by
  a simple computation along the lines of \cite{Eu2}).

  (ii) For graded degree reasons, using \eqref{hsce1},\eqref{hsce2}, and the
  fact that $\uHH^{4}(A) \cong \uHH^0(A)(-2h)$, the triples $(1,3,3)$
  and $(2,2,3)$ of multiplications \eqref{cecuppos} are zero.  Then,
  $\uHH^1(A) \cup \uHH^1(A) = 0$ since $Eu \cup Eu = 0$, by
  graded-commutativity.  The final statement then follows from Theorem
  \ref{cyfhhthm}.
\end{proof}
We now describe explicitly the Connes and BV differentials.  For this,
we fix the isomorphism
$\D: \uHH_{\bullet}(A) \cong \uHH^{3-\bullet}(A)$ of Theorem
\ref{cyfhhthm}, and use the elements
$Eu \in \uHH^{4m+1}(A), Eu^* \in \uHH^{4m+2}(A)$, and $1^* \in \uHH^{4m+3}(A)$.  Here,
the notation is a bit abusive, since really
$Eu^* \in \uHH^6(A), 1^* \in \uHH^7(A)$ using Theorem \ref{cyfhhthm},
but we identify these elements with their images under the
periodicity.  We describe all elements of $\uHH^{4m+2}(A), \uHH^{4m+3}(A)$ by
$Eu^*/z, 1^*/z$ for $z \in \uHH^0(A)$, which refers to the unique
elements so that $z \cup Eu^*/z = Eu^*$ and $z \cup 1^*/z = 1^*$.
\begin{thm}
With the above identifications, the BV differential is given by, for all
$m \in \Z$,
\begin{gather}
\Delta_{2m} = 0, \\
\Delta_{4m+1}(z Eu) = (\deg(z)+4 - 2hm) z \cdot (1')^{\cup m},  
\\ \Delta_{4m+3} (1^*/z) = (2h(1-m) - 4 - \deg(z)) Eu^*/z.
\end{gather}
In particular, $\Delta$ is graded selfadjoint, i.e., $\uHH^\bullet$ is
BV Frobenius.
\end{thm}
\begin{proof}
  We use the Cartan identity \eqref{cartan} in the case $a = Eu$:
  $B i_{Eu} + i_{Eu} B = \LL_{Eu}$.  Also, it is easy to check (as in
  e.g.~\cite{Eu3}) that $\LL_{Eu}(f) = \deg(f) \cdot f$ for all
  $f \in \HH_{\bullet}(A)$ (and hence for $\uHH_\bullet(A)$ as well).
  From this (using \eqref{hsce1}, \eqref{hsce2}, the fact that the
  Calabi-Yau shift is $4$, and the vanishing of $\uHH^1(A) \cup \uHH^3(A)$) we compute $B|_{\HH_0(A)}$:
\begin{equation}
\D(Eu \cup \Delta(1^*/z)) = (B i_{Eu} + i_{Eu} B)(\D(1^*/z)) = (2h-4- \deg(z)) \cdot 1^*/z,
\end{equation}
which implies that $\Delta_3(1^*/z) = (2h-2- \deg(z) ) Eu^*/z$, using Theorem \ref{cprepthm}.(ii).  Then, $\Delta^2 = 0$ implies that $\Delta_2 = 0$. 

Inductively, we claim that $\Delta_{2m} = 0$ and $\Delta_{2m+1}$ is an
isomorphism, for all $m \leq 1$.  Assume that $\Delta_{2m} = 0$ and
$\Delta_{2m+1}$ is an isomorphism. We need only show that
$\Delta_{2m-1}$ is also an isomorphism. This follows from
$\LL_{Eu} = B_{3-2m} i_{Eu} + i_{Eu} B_{4-2m} = i_{Eu} B_{4-2m}$,
together with the fact that $\LL_{Eu}$ and $i_{Eu}$ are isomorphisms
(here, we use that the degrees of $\uHH_0(A)$ are between $0$ and
$2h-6$).  At the same time, we may deduce the desired formulas (since
$\LL_{Eu}$ multiplies by degree).

It remains only to prove that the formulas still hold in Hochschild
degrees $> 3$.  To show this, we consider the formula
\begin{equation}
\LL_{1'} = B i_{1'} - i_{1'} B,
\end{equation}
for $1' \in \uHH^4(A)$ the periodicity element.  From this, we may
compute that, applied to degrees $\leq 3$, $\LL_{1'}$ kills even
degrees, and on odd degrees acts by $\LL_{1'}(z Eu) = 2h z (1')^{m}$
(where here $|z Eu| = 4m+1$), and $\LL_{1'}(1^*/z) = 2h Eu^*/z$.
Since $\LL_{1'}$ is a derivation and $\LL_{1'}(1') = 0$, we deduce the
desired result.
\end{proof}

\section{Periodic group algebras of finite groups}
As mentioned already, for any finite group $G$, the group algebra
$\k[G]$ is Frobenius, and in fact symmetric (hence, Calabi-Yau
Frobenius). It is natural to ask when such group algebras are
periodic.

Certainly, if $\k[G]$ is periodic, then its Hochschild cohomology is
periodic.  It is well known that one has the following formula for
Hochschild cohomology, as an abstract graded $\k$-module:
\begin{equation} \label{hhkge}
\HH^\bullet(\k[G], \k[G]) = \bigoplus_{\underset{\text{with 
representative $c_i \in C_i$}}{\text{conjugacy classes } C_i}}
H^\bullet(Z_G(c_i), \k),
\end{equation}
where $Z_G(c_i)$ is the centralizer of $c_i$ in $G$, and the
$H^\bullet(H,\k)$ denotes the \textit{group cohomology} of $H$ with
coefficients in $\k$. (For an explanation, see Proposition
\ref{projprop}, where we give a refined version.)

Hence, in order for the Hochschild cohomology to be periodic, it must
be that the numbers of generators of the cohomology groups of $G$ are
bounded.  Let us now set $\k := \Z$. Then, the classical Suzuki-Zassenhaus
theorem classifies all such groups. 
These groups are  
those such that all abelian subgroups are cyclic, and they fall into
six explicit families (cf.~p.~150 of \cite{AM}). Moreover, these
all have periodic group cohomology.  Since this property
is preserved under taking subgroups, we deduce the (probably well known)
\begin{prop} The group algebra $\Z[G]$ of a finite group $G$ has
  periodic Hochschild cohomology iff all abelian subgroups of $G$ are
  cyclic.  For such groups, $\k[G]$ has periodic Hochschild cohomology
  (relative to $\k$), for all commutative rings $\k$.
\end{prop}

We would like to know if such group algebras are in fact periodic
Calabi-Yau Frobenius algebras (since they are symmetric, it is enough
to check if they are periodic Frobenius).  This is stronger than
having periodic Hochschild cohomology, since we actually need
$\Omega^n \k[G] \simeq \k[G]$ for some $n \geq 1$.  This would be
satisfied if we could show that such $\k[G]$ have periodic
resolutions.

Fortunately, there is a very similar classical result of Swan (which
also used the (mod-$p$) classification of periodic groups):
\begin{thm}\cite{Sw1} Let $R := \Z[S^{-1}]$ for some set $S$ of primes.
Let $G$ be a finite group. Then, there is a periodic resolution of $R$ as
an $R[G]$ module, i.e.,
\begin{equation}
  0 \rightarrow R \into P_{n-1} \rightarrow P_{n-2} \rightarrow 
\cdots \rightarrow P_{0} \onto R \rightarrow 0,
\end{equation}
iff $G$ has periodic group cohomology with coefficients in $R$.  
\end{thm}
We deduce the following:
\begin{thm} \label{gpalgthm}
The group algebra $\Z[G]$ is periodic Calabi-Yau Frobenius iff $G$ has
the property that all abelian subgroups are cyclic. Equivalently,
$\Z[G]$ is periodic CY Frobenius of period $n$ iff its Hochschild
cohomology (or the group cohomology of $G$) is periodic of period $n$.
\end{thm}
\begin{proof}
  Using the above results, it is enough to show how to explicitly pass
  between a projective resolution of $\k$ as a $\k[G]$-module, and a
  projective resolution of $\k[G]$ as a $\k[G]^e$-module, in such a
  way as to preserve periodicity of a given period.  More generally,
  we prove the following proposition.
\end{proof}
\begin{prop}\label{projprop}
  Let $H=(H,\mu,\Delta,\eta, \epsilon,S)$ be any Hopf algebra over
  $\k$ which is a projective $\k$-module.  Let
  $H_\Delta := (1 \o S)\circ \Delta(H) \subset H \o H^{\op}$.  Then,
  given any projective resolution
  $P_\bullet$ of $\k$ as an $H$-module,
  $\Ind_{H_\Delta}^{H^e} P_\bullet$ is a projective
  resolution of $H$ as an $H$-bimodule, which is split as a complex
of left $H$-modules.

Conversely, if $Q_\bullet$ is a resolution of $H$ as an $H$-bimodule which
is split as a complex of right $H$-modules,
then $Q_\bullet \o_{H} \k$ is a projective resolution of $\k$ as an
$H$-module.
\end{prop}
To remove the left-right asymmetry, note that a sequence of
(projective) $H$-bimodules which are split as left $H$-modules can
have the $H$-bimodule action twisted to $a \star m \star b := S(b) m S(a)$
to make them split as right $H$-modules rather than as left
$H$-modules. If we apply $\o_H \k$ to the twisted version of $H$,
we still obtain $\k$.
\begin{proof} 
  We only prove the first statement, since the last one is easy.
    The proof is essentially a refinement of the usual
  proof of \eqref{hhkge}.  We claim that (1)
  $\Ind_{H_\Delta}^{H \o H^{\op}} \k \cong H$ as $H$-bimodules, and
  (2) with the left $(H \o 1)$ and right $H_\Delta$ actions,  
  $H^e$ is isomorphic
  as an $H$-bimodule to $H \o H$ with the usual outer $H$-bimodule structure.
By part (2) of the claim, and the fact that $H$ is projective over
  $\k$, we deduce that $H \o H^\op$ is a projective $H_\Delta$ module.
  Since $\Ind_{H_\Delta}^{H \o H^\op}$ is the functor
  $H^e \o_{H_\Delta} -$, we obtain the desired result.

  To prove the claim, consider the $\k$-linear maps
\begin{gather}
\phi: H \o H \rightarrow H \o H^\op, \quad \phi(g \o h)
= g \cdot (1 \o S)\circ\Delta(h), \\
\psi: H \o H^\op \rightarrow H \o H, \quad \psi(g \o h) = g \cdot (S \otimes 1) \circ \Delta(S^{-1} h) 
\end{gather}
By the antipode identity, coassociativity, and the counit condition,
$\phi \circ \psi = \Id = \psi \circ \phi$. On the other hand, $\phi$
intertwines the right $H$-module structure on $H \o H$ with the right
$H_\Delta = H$-module structure on $H \o H^\op$, and $\psi$
intertwines in the opposite direction; also, both intertwine the
standard left $H$-module structure.  
So, we obtain part (2).  Part (1) 
 then follows by explicit (easy) computation.
\end{proof}
\begin{cor} A Hopf algebra has a periodic bimodule resolution which is
  split as a complex of right modules iff its augmentation module has
  a periodic left module resolution.  A Hopf algebra which is a
  Frobenius algebra is periodic Frobenius iff its augmentation module
  $\k$ satisfies $\Omega^n \k \simeq \k$ in the stable left-module
  category.
\end{cor}
As remarked earlier, if $\k$ is a PID, any Hopf algebra which is fg as
a $\k$-module is automatically Frobenius \cite{LS}, so we can remove
the Frobenius assumption from the corollary in this case (replacing
with fg projective over $\k$).
\begin{proof}
  The first assertion follows immediately from the proposition.  For
  the second, we show that one has a stable module equivalence
  $\ds \Omega^n \k \mathop{\simeq}^H \k$ iff one has a stable bimodule
  equivalence $\ds \Omega^n H \mathop{\simeq}^{H^e} H$. We can use
  the functors $- \uo_H \k, H^e \uo_{H_\Delta} -$ to
  achieve this. 
\end{proof}
The periodic groups described by the theorem include
all groups which act freely on spheres:
\begin{thm} \label{fingpthm} Let $G$ be any finite group which acts
  freely and orientation-preserving on a sphere $S^\ell$ with $\ell \geq 1$.
Then, for some $r \geq 1$, the
  group algebra $\k[G]$ is a periodic Calabi-Yau Frobenius algebra, of
  dimension $\frac{\ell+1}{r} - 1$ and period $\frac{\ell+1}{r}$.
\end{thm}
This theorem follows from e.g.~Lemma 6.2 of \cite{AM} (a spectral
sequence argument showing that group cohomology is periodic in this
case), by using Theorem \ref{gpalgthm}. We give, however, a simple
topological proof in the Appendix, in the case when $G$ acts
cellularly on a finite CW complex homeomorphic to $S^\ell$.
\begin{cor}
  For any finite subgroup $G < SO(\ell+1) := SO(\ell+1,\R)$ for any
  $\ell \geq 1$, $\k[G]$ is periodic Calabi-Yau Frobenius.
\end{cor}
We note that this corollary may also be deduced from the version
of the theorem proved in the appendix, where $G$ acts cellularly.
To do this, we form a CW decomposition of $S^{\ell+1}$ by
geodesic codimension-one slices, fixed under the orbit of $G$, which
separate a given point $x$ from all of its orbits under $G$.

\appendix
\section{Finite groups acting freely on spheres} \label{topapp}
In this section, we provide an elementary proof of Theorem
\ref{fingpthm} in the case that $G$ acts cellularly on a finite
CW complex homeomorphic to $S^\ell$.
 Our proof avoids the use of the classification of periodic groups,
 and is purely topological.

The main idea of the proof is to construct a bimodule resolution of
$\k[G]$ by constructing a CW complex which is homotopic to $G$, and
which admits a free action of $G \times G^{\op}$, such that the
induced $G \times G^{\op}$-module structure on the cellular homology, $\k[G]$,
is the standard bimodule structure.  The CW complex will have
finitely many cells of each dimension, and the resulting cellular
chain complex will be periodic.  Thus, this complex yields a periodic
bimodule resolution of $\k[G]$.

To do this, we need the following simple topological lemma:
\begin{ntn} Let $D^n, S^n$ denote the $n$-dimensional disc and sphere,
respectively.
\end{ntn}
\begin{lemma} \label{fglem}
  Let $m, n \geq 0$ be any integers.  Consider the topological space
  $X := D^{m+1} \times S^n$, and let $f: \partial X \rightarrow S^m$
  be the attaching map which is the first projection of
  $S^{m} \times S^n$ to $S^{m}$.  Then, the glued topological
  space $X \cup_{f} S^m$ is homeomorphic to $S^{n+m+1}$.
\end{lemma}
In the special case that $n = 0$, the above construction is the standard
way to build $S^{m+1}$ out of $S^m$: we attach two hemispheres $D^{m+1}$ to
$S^m$ placed at the equator.
\begin{proof}
Let us view $S^{n+m+1} \subseteq \R^{n+m+2}$ as the unit sphere. Let
$Y^\circ \subset S^{n+m+1}$ be the subset 
\begin{equation}
Y^\circ := \{(x_0, \ldots, x_{n+m+1}) \in S^{n+m+1} \subset \R^{n+m+2}: \sum_{i=0}^{m} |x_i|^2 < 1\}.
\end{equation}
Now, $Y^\circ$ is an open subset of $S^{n+m+1}$ homeomorphic to
$D^{m+1} \times S^n$, under the map 
\begin{equation}
D^{m+1} \times S^n \ni (\mathbf{x},\mathbf{y}) \mapsto
(\mathbf{x}, \frac{1}{\sqrt{1- \norm x^2}} \mathbf{y}).
\end{equation}
The complement of $Y^\circ$ is the subset $S^{m} \times \{0\} \subset
S^{n+m+1}$.  That is, $\partial(Y^\circ) \cong S^m$, and the attaching
map $S^{m+1} \times S^n \rightarrow \partial(Y^\circ)$ is $f$.
\end{proof}
In fact, we will not only use the statement of the lemma, but (a
cellular version of) the explicit homeomorphism given in the proof.
\begin{proof}[Proof of the cellular version of Theorem \ref{fingpthm}]
  We will construct the topological space
  $X := S^\infty \times_{S^\infty/G} S^\infty$ as an explicit CW
  complex with finite $n$-skeleta, whose associated cellular chain
  complex is $\ell+1$-periodic. Moreover, the group $G \times G$ will act
  freely and cellularly (since $G$ is finite, this
  means that the group sends each $d$-dimensional cell onto another
  $d$-dimensional cell).  This gives the desired result by the remarks
  at the beginning of the subsection.

  Let $Y$ denote the finite CW complex with $Y \cong S^{\ell}$ that we
  are given.  First, we construct from this a topological space
  $Z \cong D^{\ell+1}$ that extends the action, viewing $Y$ as the
  boundary of $Z$.  To do this, let
  $Z := (Y \times [0,1]) / (Y \times \{0\})$, and let $G$ act in the
  obvious way (preserving the second component).  We will think of $Z$
  as the solid unit disc in $\R^{\ell+1}$ and of $Y$ as its boundary.
  We will also view $Z$ as a mere topological space isomorphic to
  $D^{\ell+1}$, i.e., a single $\ell+1$-cell with a $G$-action, for
  the purpose of constructing complexes.

  Now, set $W_\ell := Y$.  We inductively construct a copy $W$ of $S^\infty$ 
by attaching (viewing $Z^k$ on the LHS as a single $k(\ell+1)$-cell):
\begin{gather} \label{attmaps}
Z \times Y \cong D^{\ell+1} \times Y \tra^{att.} Y = W_\ell, \\
Z^2 \times Y \cong D^{2(\ell+1)} \times Y \tra^{att.} 
W_{(\ell+1)+\ell}, \\
Z^3 \times Y \cong D^{3(\ell+1)} \times Y \tra^{att.} 
W_{2(\ell+1)+\ell}, \ldots,
\end{gather}
where ``att.'' means an attaching map (so NOT a map of topological
spaces). We
define these attaching maps to be the composition of the first
projection $S^{k(\ell+1)-1} \times Y \onto S^{k(\ell+1)-1}$ with the
homeomorphism
$S^{k(\ell+1)-1} \iso W_{(k-1)(\ell+1)+\ell} = W_{k(\ell+1)-1}$, which
exists inductively by Lemma \ref{fglem}.  Constructed this way, each
homeomorphism $S^{k(\ell+1)+\ell} \iso W_{k(\ell+1)+\ell}$ has image
in the same sum of top cells of $Z^k \times Y$: it is the sum of the
$\ell$ cells of $Y$ which make up $S^{\ell} \iso Y$.

Thus, on the level of chains, if we label the cells of $W_{\infty}$ in
degree $k(\ell+1)+p$ by $c_{k,p} = Z^k \times c_{0,p}$, where
$c_{0,p}$ are the cells of $Y$ for $0 \leq p \leq \ell$, the complex
$C_\bullet(W_\infty)$ is a periodic free complex.

By construction, $W_{\infty}$ has a free action of $G$.  Now, finally,
set $X := W_{\infty} \times G$, and let us view $X$ as the
homeomorphic space $W_{\infty} \times_{W_{\infty}/G} W_{\infty}$, and
equip $X$ with the resulting free cellular action of $G \times G$.  We
then have that $X$ is a topological space whose homology is $\k[G]$
with the usual bimodule action.  We deduce that $C_\bullet(X)$ is a
free periodic resolution of $\k[G]$.
\end{proof}
\section{Frobenius algebras over general commutative base rings}
\label{gfs} In this subsection, we will extend some results known for
Frobenius algebras over fields to a relative context, using
e.g.~\cite{ARS} as a reference for the usual results.  

We first deduce a relative selfinjectivity for $A$.
\begin{ntn} For any $\k$-algebra $A$ which is fg
  projective as a $\k$-module, let  \\ $^*$ denote the functor
  $*: \Hom_\k(-,\k): A-mod \rightarrow A^{\op}-mod$.
\end{ntn}
\begin{defn}
  Call an $A$-module $I$ which is fg projective over
  $\k$ ``injective relative to $\k$'' or ``relatively injective'' if,
  for any $A$-modules $M$, $N$ which are projective over $\k$, and any
  injection $M \into N$ whose cokernel is a projective $\k$-module
  (i.e., the injection is $\k$-split), there exists a unique dotted
  arrow completing any diagram of the following form:
\begin{equation}
\xymatrix{
I & \\
M \ar@{^{ (}->}[r] \ar[u] & \ar@{-->}[ul] N 
}
\end{equation}
\end{defn}
The following helps explain the meaning of relative injectivity:
\begin{prop} \label{riap} Let $A$ be any algebra over a commutative ring $\k$.
If $M$ is any module over $A$ which is fg projective as a $\k$-module,
then any relatively injective $A$-module is
acyclic for the functor $\Hom_A(M, -)$. (I.e., such $I$ satisfy $\Ext_A^i(M, I) = 0$ for all $i \geq 1$.) 
\end{prop}
\begin{proof}
Let us take a projective resolution $P_\bullet$ of $M$, i.e.,
$\cdots \rightarrow P_1 \rightarrow P_0 \onto M$.  Since $M$ is projective over $\k$, the resolution is $\k$-split.  Thus, the relative injectivity property
will guarantee that there is no positive homology of $\Hom(P_\bullet, I)$.
\end{proof}

The main use of ``relative to $\k$'' is to make the dualization
$*: M \mapsto \Hom_\k(M, \k)$ behave well:
\begin{lemma} Assume $A$ is a $\k$-algebra which is a
  fg projective $\k$-module. Then, 
\begin{enumerate}
\item[(i)] The dualization $*: A-mod \rightarrow A^{\op}-mod$ is a
  contravariant functor which sends projective modules to
  modules which are injective relative to $\k$, and vice-versa.
\item[(ii)] There is a functorial isomorphism $M^{**} \cong M$, for
  $M$ fg projective over $\k$.  That is, the
  restriction of $*$ to fg $A$-modules which are
  projective over $\k$ (in both the domain and codomain) is an
  anti-equivalence of categories, and $* \circ * \simeq \Id$.
\end{enumerate}
\end{lemma}
The proof is just as in the case where $\k$ is a field, so we omit it.
\begin{cor}\label{eric}
  If $A$ is as in the lemma, then $A-mod$ has enough relatively
  injectives in the following sense: for every fg $A$-module $M$ which is
  projective over $\k$, there exists a relatively injective module $I$
  and a $\k$-split injection $M \into I$.
\end{cor}
\begin{proof}
  For any $\k$-projective $M \in A-mod$, pick a surjection $P \onto M^*$ in the
  category $A^{\op}-mod$.  This is $\k$-split. Then, dualizing, one
  obtains a $\k$-split injection $M \into P^*$, and $P^*$ is
  relatively injective.
\end{proof}
As a corollary, we also deduce the relative selfinjectivity for Frobenius algebras:
\begin{cor} \label{relinjc} If $A$ is a Frobenius algebra over $\k$,
  then $A$ is relatively injective as an $A$-module.  Moreover, all
  projectives are relatively injective and vice-versa.
\end{cor}
\begin{proof}
We know that $A^*$ is isomorphic to $A$ as an $A$-module (using the invariance
and nondegeneracy of the pairing). Now, as an $A^{\op}$-module, $A$ is projective;
hence $A^*$ is relatively injective as an $A$-module.

For the last statement, we use the fact that all projectives are
direct summands of free modules, and hence are direct summands of 
relatively injective modules, and are hence relatively injective.  For
the converse,  for any relatively injective module $M$,
 $M^*$ is a projective
$A^{\op}$-module, hence a relatively injective $A^{\op}$-module, and
hence a projective $A$-module.
\end{proof}
Next we show that the other duality $\vee$ also behaves well for
Frobenius algebras in the relative context.  All of the following
results are more generally true for \emph{relatively selfinjective}
algebras, which we define as algebras $A$ satisfying the conclusion of
Corollary \ref{eric}: they are fg projective over $\k$ and relatively
injective (equivalently, $A^*$ is projective, i.e., all projectives
are relatively injective and vice-versa). The same proofs apply.
\begin{prop} \label{frobdualprop}
Suppose $A$ is Frobenius over $\k$. Then,
\begin{enumerate}
\item[(i)] The functor $\vee$ restricts to a functor on
full subcategories:
\begin{equation}
\{\text{fg $A$-modules which are projective over $\k$}\}
\leftrightarrow \{\text{fg $A^\op$-modules which are projective over $\k$}\},
\end{equation}
and $\vee \circ \vee \simeq \Id$ on these subcategories.
\item[(ii)] The functor $\vee$ preserves \textbf{exact} 
$\k$-split complexes of fg projective $\k$-modules.
\end{enumerate}
\end{prop}
\begin{proof}[Proof of Proposition \ref{frobdualprop}]
(i)
First, let us show that, if $M \in A-mod$ is fg
projective over $\k$, then so is $M^\vee$:
\begin{multline} \label{sfe1} M^\vee = \Hom_A(M, A) \cong \Hom_A(M,
  (A^*)^*) = \Hom_A(M, \Hom_\k(A^*, \k)) \\ \mathop{\cong}^{adj.}
  \Hom_\k(A^* \o_A M, \k) = (A^* \o_A M)^* \cong (_{\eta^{-1}} M)^*,
\end{multline}
which is f.g. projective over $\k$. Moreover, this is functorial in
$M$, and we deduce that $\vee \circ \vee \simeq \Id$.

(ii) For bounded-below complexes, this follows from the fact
(Proposition \ref{riap}) that $\Ext^i_A(M, A) = 0$ for all $i \geq 1$.
Then, for any unbounded $\k$-split exact complex of projectives,
we truncate at an arbitrary place.
\end{proof}
\begin{cor} \label{smdc}
Let $A$ be a Frobenius algebra over $\k$.
We have mutually inverse 
autoequivalences $\Omega, \Omega^{-1} : \uAmk \rightarrow \uAmk$, 
which yield exact sequences
\begin{equation} \label{omms}
0 \rightarrow \Omega M \into P \onto M \rightarrow 0, 
\end{equation}
with $P$ a fg projective $A$ module,
for all fg $A$-modules $M$ which are projective over $\k$.
\end{cor}
\begin{proof}
Choose, for every module $M$, a sequence \eqref{omms}, and
  similarly a sequence of the form
\begin{equation}
  0 \rightarrow M \into I \onto \Omega^{-1}(M) \rightarrow 0,
\end{equation}
with $I$ satisfying the relative injectivity property.  (To obtain
such a sequence, form a sequence of the form \eqref{omms} for $M^*$ in
the category $A^{\op}-mod$, and then dualize.)  Then, using the fact
that a map $M \rightarrow N$ factors through a specific injection
$M \into I$ for $I$ relatively injective iff it factors through any
other injection into a relatively injective (both are true iff $M$
factors through all injections), it is straightforward to finish using
the same arguments as in the case when $\k$ is a field (see, e.g.,
\cite{ARS}, \cite{ASS}).
\end{proof}
\begin{cor} \label{extcor}
Let $A$ be any Frobenius algebra over $\k$. For any degree $i \geq 1$, and any 
fg $A$-modules
$M, N$ which are projective over $\k$,
\begin{equation}
\Ext^i_A(M, N) \cong \uHom_A(\Omega^i M, N).
\end{equation}
\end{cor}
\begin{proof}
$\Ext^i_A(M, N)$ can be computed using a projective resolution of $M$,
\begin{equation}
P_{i+1} \tra^{d_{i+1}} P_i \tra^{d_i} \cdots \tra^{d_1} P_0 \onto M,
\end{equation}
by taking the quotient
\begin{equation}
\{f \in \Hom_{A^e}(P_{i}, N) \mid f \circ d_{i+1} = 0\}/ \{g \circ d_{i} 
\}_{g \in \Hom_{A^e}(P_{i-1}, N)}.
\end{equation}
Now, factor $P_i \rightarrow P_{i-1}$ as
$P_i \onto \Omega^i M \into P_{i-1}$ (for
$\Omega^i M = \imm d_i \cong \coker d_{i+1}$), so that we may now
write $\Ext_A^i(M, N)$ as
$\Hom_{A}(\Omega^i M, N)/\{\text{morphisms factoring through
  $\Omega^i M \into P_{i-1}$}\}$.
Then, by the observation from the proof of Corollary \eqref{smdc} that
morphisms factor through $\Omega^i M \into P_{i-1}$ iff they factor
through any other injection to an injective $A$-bimodule relative to $\k$
(and injective $A$-bimodules relative to $\k$ are the same as
fg projective $A$-bimodules by Corollary \ref{relinjc}), we obtain the
desired result.
\end{proof}
\begin{rem} When $\k$ is a field, the above actually endows the stable
  module category with the structure of a triangulated category, which
  is the quotient of the derived category of fg $A$-modules by finite
  complexes of projective $A$-modules.  However, this is NOT true for
  general $\k$ (we had to restrict to the non-abelian subcategory of
  $\k$-projectives before doing anything).
\end{rem}

\bibliographystyle{amsalpha}
{\small

\bibliography{cyfrob}

}
\end{document}